\documentclass[11pt,twoside,a4paper]{article}
\usepackage{amssymb}
\usepackage[english]{babel}
\newcommand{\tun}{\begin{picture}(5,0)(-2,-1)
\put(0,0){\circle*{2}}
\end{picture}}
\newcommand{\tdeux}{\begin{picture}(7,7)(0,-1)
\put(3,0){\circle*{2}}
\put(3,0){\line(0,1){5}}
\put(3,5){\circle*{2}}
\end{picture}}
\newcommand{\ttroisdeux}{\begin{picture}(5,12)(-2,-1)
\put(0,0){\circle*{2}}
\put(0,0){\line(0,1){5}}
\put(0,5){\circle*{2}}
\put(0,5){\line(0,1){5}}
\put(0,10){\circle*{2}}
\end{picture}}
\newcommand{\tquatrecinq}{\begin{picture}(5,19)(-2,-1)
\put(0,0){\circle*{2}}
\put(0,0){\line(0,1){5}}
\put(0,5){\circle*{2}}
\put(0,5){\line(0,1){5}}
\put(0,10){\circle*{2}}
\put(0,10){\line(0,1){5}}
\put(0,15){\circle*{2}}
\end{picture}}
\newcommand{\tcinqquatorze}{\begin{picture}(5,26)(-2,-1)
\put(0,0){\circle*{2}}
\put(0,0){\line(0,1){5}}
\put(0,5){\circle*{2}}
\put(0,5){\line(0,1){5}}
\put(0,10){\circle*{2}}
\put(0,10){\line(0,1){5}}
\put(0,15){\circle*{2}}
\put(0,15){\line(0,1){5}}
\put(0,20){\circle*{2}}
\end{picture}}
\newcommand{\ttroisun}{\begin{picture}(15,8)(-5,-1)
\put(3,0){\circle*{2}}
\put(-1,0){$\vee$}
\put(6,7){\circle*{2}}
\put(0,7){\circle*{2}}
\end{picture}}
\newcommand{\tquatreun}{\begin{picture}(15,12)(-5,-1)
\put(3,0){\circle*{2}}
\put(-1,0){$\vee$}
\put(6,7){\circle*{2}}
\put(0,7){\circle*{2}}
\put(3,7){\circle*{2}}
\put(3,0){\line(0,1){7}}
\end{picture}}
\newcommand{\tquatredeux}{\begin{picture}(15,15)(-5,-1)
\put(3,0){\circle*{2}}
\put(-1,0){$\vee$}
\put(6,7){\circle*{2}}
\put(0,7){\circle*{2}}
\put(0,14){\circle*{2}}
\put(0,7){\line(0,1){7}}
\end{picture}}
\newcommand{\tquatretrois}{\begin{picture}(15,15)(-5,-1)
\put(3,0){\circle*{2}}
\put(-1,0){$\vee$}
\put(6,7){\circle*{2}}
\put(0,7){\circle*{2}}
\put(6,14){\circle*{2}}
\put(6,7){\line(0,1){7}}
\end{picture}}
\newcommand{\tquatrequatre}{\begin{picture}(15,14)(-5,-1)
\put(3,5){\circle*{2}}
\put(-1,5){$\vee$}
\put(6,12){\circle*{2}}
\put(0,12){\circle*{2}}
\put(3,0){\circle*{2}}
\put(3,0){\line(0,1){5}}
\end{picture}}
\newcommand{\tcinqun}{\begin{picture}(20,8)(-5,-1)
\put(3,0){\circle*{2}}
\put(-1,0){$\vee$}
\put(6,7){\circle*{2}}
\put(0,7){\circle*{2}}
\put(3,0){\line(2,1){10}}
\put(3,0){\line(-2,1){10}}
\put(-7,5){\circle*{2}}
\put(13,5){\circle*{2}}
\end{picture}}
\newcommand{\tcinqdeux}{\begin{picture}(15,14)(-5,-1)
\put(3,0){\circle*{2}}
\put(-1,0){$\vee$}
\put(6,7){\circle*{2}}
\put(0,7){\circle*{2}}
\put(3,7){\circle*{2}}
\put(3,0){\line(0,1){7}}
\put(0,7){\line(0,1){7}}
\put(0,14){\circle*{2}}
\end{picture}}
\newcommand{\tcinqtrois}{\begin{picture}(15,15)(-5,-1)
\put(3,0){\circle*{2}}
\put(-1,0){$\vee$}
\put(6,7){\circle*{2}}
\put(0,7){\circle*{2}}
\put(3,7){\circle*{2}}
\put(3,0){\line(0,1){7}}
\put(3,7){\line(0,1){7}}
\put(3,14){\circle*{2}}
\end{picture}}
\newcommand{\tcinqquatre}{\begin{picture}(15,14)(-5,-1)
\put(3,0){\circle*{2}}
\put(-1,0){$\vee$}
\put(6,7){\circle*{2}}
\put(0,7){\circle*{2}}
\put(3,7){\circle*{2}}
\put(3,0){\line(0,1){7}}
\put(6,7){\line(0,1){7}}
\put(6,14){\circle*{2}}
\end{picture}}
\newcommand{\tcinqsix}{\begin{picture}(15,15)(-7,-1)
\put(3,0){\circle*{2}}
\put(-1,0){$\vee$}
\put(6,7){\circle*{2}}
\put(0,7){\circle*{2}}
\put(-4,7){$\vee$}
\put(3,14){\circle*{2}}
\put(-3,14){\circle*{2}}
\end{picture}}
\newcommand{\tcinqsept}{\begin{picture}(15,8)(-5,-1)
\put(3,0){\circle*{2}}
\put(-1,0){$\vee$}
\put(6,7){\circle*{2}}
\put(0,7){\circle*{2}}
\put(2,7){$\vee$}
\put(3,14){\circle*{2}}
\put(9,14){\circle*{2}}
\end{picture}}
\newcommand{\tcinqhuit}{\begin{picture}(15,26)(-5,-1)
\put(3,0){\circle*{2}}
\put(-1,0){$\vee$}
\put(6,7){\circle*{2}}
\put(0,7){\circle*{2}}
\put(0,14){\circle*{2}}
\put(0,7){\line(0,1){7}}
\put(0,21){\circle*{2}}
\put(0,14){\line(0,1){7}}
\end{picture}}
\newcommand{\tcinqneuf}{\begin{picture}(15,26)(-5,-1)
\put(3,0){\circle*{2}}
\put(-1,0){$\vee$}
\put(6,7){\circle*{2}}
\put(0,7){\circle*{2}}
\put(6,14){\circle*{2}}
\put(6,7){\line(0,1){7}}
\put(6,21){\circle*{2}}
\put(6,14){\line(0,1){7}}
\end{picture}}
\newcommand{\tcinqcinq}{\begin{picture}(15,19)(-5,-1)
\put(3,0){\circle*{2}}
\put(-1,0){$\vee$}
\put(6,7){\circle*{2}}
\put(0,7){\circle*{2}}
\put(6,14){\circle*{2}}
\put(6,7){\line(0,1){7}}
\put(0,14){\circle*{2}}
\put(0,7){\line(0,1){7}}
\end{picture}}
\newcommand{\tcinqdix}{\begin{picture}(15,19)(-5,-1)
\put(3,5){\circle*{2}}
\put(-1,5){$\vee$}
\put(6,12){\circle*{2}}
\put(0,12){\circle*{2}}
\put(3,0){\circle*{2}}
\put(3,0){\line(0,1){5}}
\put(3,7){\line(0,1){5}}
\put(3,12){\circle*{2}}
\end{picture}}
\newcommand{\tcinqonze}{\begin{picture}(15,26)(-5,-1)
\put(3,5){\circle*{2}}
\put(-1,5){$\vee$}
\put(6,12){\circle*{2}}
\put(0,12){\circle*{2}}
\put(3,0){\circle*{2}}
\put(3,0){\line(0,1){5}}
\put(0,12){\line(0,1){7}}
\put(0,19){\circle*{2}}
\end{picture}}
\newcommand{\tcinqdouze}{\begin{picture}(15,26)(-5,-1)
\put(3,5){\circle*{2}}
\put(-1,5){$\vee$}
\put(6,12){\circle*{2}}
\put(0,12){\circle*{2}}
\put(3,0){\circle*{2}}
\put(3,0){\line(0,1){5}}
\put(6,12){\line(0,1){7}}
\put(6,19){\circle*{2}}
\end{picture}}

\newcommand{\tcinqtreize}{\begin{picture}(5,26)(-2,-1)
\put(0,0){\circle*{2}}
\put(0,0){\line(0,1){7}}
\put(0,7){\circle*{2}}
\put(0,7){\line(0,1){7}}
\put(0,14){\circle*{2}}
\put(-4,14){$\vee$}
\put(-3,21){\circle*{2}}
\put(3,21){\circle*{2}}
\end{picture}}

\newcommand{\tdun}[1]{\begin{picture}(10,0)(-2,-1)
\put(0,0){\circle*{2}}
\put(3,-2){\tiny #1}
\end{picture}}
\newcommand{\tddeux}[2]{\begin{picture}(12,5)(0,-1)
\put(3,0){\circle*{2}}
\put(3,0){\line(0,1){5}}
\put(3,5){\circle*{2}}
\put(6,-2){\tiny #1}
\put(6,3){\tiny #2}
\end{picture}}
\newcommand{\tdtroisun}[3]{\begin{picture}(20,8)(-5,-1)
\put(3,0){\circle*{2}}
\put(-1,0){$\vee$}
\put(6,7){\circle*{2}}
\put(0,7){\circle*{2}}
\put(5,-2){\tiny #1}
\put(9,5){\tiny #2}
\put(-5,5){\tiny #3}
\end{picture}}
\newcommand{\tdtroisdeux}[3]{\begin{picture}(12,12)(-2,-1)
\put(0,0){\circle*{2}}
\put(0,0){\line(0,1){5}}
\put(0,5){\circle*{2}}
\put(0,5){\line(0,1){5}}
\put(0,10){\circle*{2}}
\put(3,-2){\tiny #1}
\put(3,3){\tiny #2}
\put(3,8){\tiny #3}
\end{picture}}

\newcommand{\tdquatrequatre}[4]{\begin{picture}(20,14)(-5,-1)
\put(3,5){\circle*{2}}
\put(-1,5){$\vee$}
\put(6,12){\circle*{2}}
\put(0,12){\circle*{2}}
\put(3,0){\circle*{2}}
\put(3,0){\line(0,1){5}}
\put(6,-3){\tiny #1}
\put(6,4){\tiny #2}
\put(9,12){\tiny #3}
\put(-5,12){\tiny #4}
\end{picture}}

\newcommand{\malv}{{\bf FQSym}}

\renewcommand{\P}{{\cal P}}
\renewcommand{\L}{{\cal L}}
\newcommand{\F}{\mathbb{F}}
\newcommand{\FF}{{\bf F}}
\newcommand{\A}{{\cal A}}
\newcommand{\h}{{\cal H}}
\newcommand{\D}{{\cal D}}
\renewcommand{\S}{\mathbb{S}}
\newcommand{\T}{\mathbb{T}}
\newcommand{\V}{{\bf V}}
\newcommand{\tdelta}{\tilde{\Delta}}
\newcommand{\totimes}{\overline{\otimes}}

\title{Bidendriform bialgebras, trees, and free quasi-symmetric functions}
\date{}
\author{L. Foissy \\
\\
{\small{\it Laboratoire de Math\'ematiques - UMR6056, Universit\'e de Reims}}\\
\small{{\it Moulin de la Housse - BP 1039 - 51687 REIMS Cedex 2, France}}
\footnote{e-mail : loic.foissy@univ-reims.fr}}
\newtheorem{lemme}{\indent Lemma}
\newtheorem{defi}[lemme]{\indent Definition}
\newtheorem{prop}[lemme]{\indent Proposition}
\newtheorem{theo}[lemme]{\indent Theorem}
\newtheorem{cor}[lemme]{\indent Corollary}

\begin{document}

\maketitle

ABSTRACT: we introduce bidendriform bialgebras,
which are bialgebras such that both product and coproduct can be split
into two parts satisfying good compatibilities. For example, 
the Malvenuto-Reutenauer Hopf algebra and the non-commutative Connes-Kreimer
Hopf algebras of planar decorated rooted trees are bidendriform bialgebras.
We prove that all connected bidendriform bialgebras are generated by their 
primitive elements as a dendriform algebra (bidendriform Milnor-Moore theorem)
and then is isomorphic to a Connes-Kreimer Hopf algebra. As a corollary, the Hopf algebra of
Malvenuto-Reutenauer is isomorphic to the Connes-kreimer Hopf algebra of planar rooted trees 
decorated by a certain set. We deduce that the Lie algebra of its primitive elements 
is free in characteristic zero (G. Duchamp, F. Hivert and J.-Y. Thibon conjecture).\\

RESUME : nous introduisons les big\`ebres bidendriformes, qui sont des big\`ebres dont
le produit et le coproduit peuvent \^etre scind\'es en deux avec de bonnes compatibilit\'es.
Par exemple, l'alg\`ebre de Hopf de Malvenuto-Reutenauer et les alg\`ebres de Hopf non-commutative
de Connes-Kreimer sur les arbres plans enracin\'es d\'ecor\'es sont des big\`ebres bidendriformes.
Nous montrons que toute big\`ebre bidendriforme connexe est engendr\'ee par ses \'el\'ements totalement
primitifs comme alg\`ebre dendriforme (version bidendriforme du th\'eor\`eme de Milnor-Moore)
et qu'elle est alors isomorphe \`a une alg\`ebre de Hopf de Connes-Kreimer.
En cons\'equence, l'alg\`ebre de Hopf de Malvenuto-Reutenauer est isomorphe \`a l'alg\`ebre
de Connes-Kreimer des arbres plans enracin\'es d\'ecor\'es par un certain ensemble. 
On en d\'eduit que l'alg\`ebre de Lie de ses \'el\'ements primitifs est libre en
caract\'eristique z\'ero (conjecture de G. Duchamp, F. Hivert et J.-Y. Thibon).\\ 

\tableofcontents

\section*{Introduction}
 
The Hopf algebra $\malv$ of Malvenuto-Reutenauer, also called Hopf algebra of free quasi-symmetric functions
(\cite{Duchamp,Loday,Loday3,Malvenuto}) has certain interesting properties. For example,
it is known that it free as an algebra and cofree as a coalgebra; it has a non-degenerate Hopf pairing;
it can be given a structure of dendriform algebra. The Hopf algebras $\h^\D$ of planar rooted trees,
introduced in \cite{Foissy2,Foissy3,Holtkamp} have a lot of similar properties: they are free as algebras
and cofree as coalgebras; they have a non-degenerate Hopf pairing (although not so explicit as the
Malvenuto-Reutenauer algebra's); they also are dendriform algebras. 
So a natural question is: are these two objects isomorphic? More precisely, is there a set $\D$ of decorations
such that $\h^\D$ is isomorphic to $\malv$?

In order to answer positively this question, we study more in details dendriform algebras.
This notion is introduced by Loday and Ronco in \cite{Loday2} and 
is studied in \cite{Aguiar,Loday,Loday3,Ronco2,Ronco}. A dendriform algebra $A$ 
is a (non unitary) associative algebra, such that the product can be split into two parts $\prec$ and $\succ$
(left and right products),with good compatibilities 
(which mean that $(A,\prec,\succ)$ is a bimodule over itself).
The notion of dendriform bialgebra or Hopf dendriform algebra is also given. These are dendriform algebra
with a coassociative coproduct $\tdelta$, satisfying good relations with $\prec$ and $\succ$.
We introduce here the notion of bidendriform bialgebra
(section \ref{sect2}): a bidendriform bialgebra $A$ is a dendriform bialgebra
such that the coproduct $\tdelta$ can be split into two parts $\Delta_\prec$ and $\Delta_\succ$
(left and right coproducts), such that $(A,\Delta_\succ,\Delta_\prec)$ is a bimodule over itself.
There are also compatibilities between the left and right coproducts and the left and right products.
As this set of axioms is self-dual, the dual of a finite-dimensional bidendriform bialgebra
is also a bidendriform bialgebra. An example of bidendriform bialgebra is $\malv$ (section \ref{sect4}),
or, more precisely, its augmentation ideal:
as $\malv$ is a dendriform algebra and it is self dual as a Hopf algebra, we can also split the coproduct
into two parts. Fortunately, the left and right coproducts defined in this way satisfy the wanted 
compatibilities with the left and right products.

We would like to give $\h^\D$ a structure of bidendriform algebra too. The method used for $\malv$
fails here: the left and right coproducts defined by duality, denoted here by $\Delta'_\prec$
and $\Delta'_\succ$, do not satisfy the compatibilities with the left and right products.
Hence, we have to proceed in a different way. For this, we consider the category of dendriform algebras
and give it a tensor product $\totimes$ (section \ref{sect3}).
This tensor product is a little bit different of the usual one,
as dendriform algebras are not unitary objects: we have to add a copy of both algebras to their tensor
products. We also define the notion of dendriform module over a dendriform algebra.
Then the notions of dendriform bialgebra and bidendriform bialgebra become more clear:
the coproduct of a dendriform bialgebra $A$ has to be a morphism of dendriform algebras from
$A$ to $A\totimes A$; the left and right coproducts of a bidendriform bialgebra $A$ have to be
morphism of dendriform modules from $A$ to $A\totimes A$.
Now, as the dendriform algebra $\h^\D$ is freely generated by the elements $\tdun{$d$}$ as a dendriform
algebra, it is possible to define a unique structure of bidendriform bialgebra over $\h^\D$
(or, more exactly, on its augmentation ideal $\A^\D$) by $\Delta_\prec(\tdun{$d$})=\Delta_\succ(\tdun{$d$})=0$
(theorem \ref{theoADdend}).

Let us now study more precisely the notion of bidendriform bialgebra. 
For a given bidendriform bialgebra $A$, we consider totally primitive elements of $A$,
that is to say elements which vanish under both $\Delta_\prec$ and $\Delta_\succ$.
We say that $A$ is connected if, for every element $a$, the various iterated coproducts all vanish
on $a$ for a great enough rank. Then, if $A$ is connected, we prove that $A$ is generated
as a dendriform algebra by its totally primitive elements (theorem \ref{theo17}).
This theorem can be seen as a bidendriform  version of the Milnor-Moore theorem (\cite{Milnor}), 
which says that a cocommutative, connected Hopf algebra is generated by its primitive elements 
(in characteristic  zero).

To precise this result, we consider the bidendriform bialgebra $\A^\D$. We prove that
its space of totally primitive elements is reduced to the space of its generators: the elements
$\tdun{$d$}$. In other terms, the triple $(Dend,coDend,Vect)$ is a good triple of operads,
with the language of \cite{Loday4}. This implies that every connected bidendriform bialgebra 
is freely generated by its totally primitive elements, so is isomorphic to a $\h^\D$ for a
well chosen $\D$ (theorem \ref{theo39}).
We apply this result to $\malv$: then, for a certain $\D$, $\malv$ is isomorphic to $\h^\D$
as a bidendriform bialgebra, and hence as a Hopf algebra. This allows us to answer a conjecture
of \cite{Duchamp}: if the characteristic of the base field is zero, then the Lie algebra of 
primitive elements of $\malv$ is free (corollary \ref{cor40}), as we already proved this result 
for $\h^\D$ in \cite{Foissy3}.

As a connected bidendriform bialgebra is isomorphic to a certain $\h^\D$, it is self-dual.
We define the notion of bidendriform pairing; for example, the pairing of \cite{Duchamp}
on $\malv$ is a bidendriform pairing, whereas the pairing of  \cite{Foissy2} on $\h^\D$ is not.
We construct another pairing on $\h^\D$, which is non-degenerate, symmetric and bidendriform.
This pairing can be studied in the same way as in \cite{Foissy2}: we give a inductive way to compute it
and a combinatorial interpretation of this pairing.
\\

{\bf Thanks.} I am grateful to Ralf Holtkamp for suggestions which greatly improve section \ref{sect3}.
\\

{\bf Notation.} $K$ is a commutative field of any characteristic.

\section{Dendriform and codendriform bialgebras}

\subsection{Dendriform algebras and coalgebras}

\begin{defi}
\textnormal{
(See \cite{Aguiar,Loday2,Loday,Loday3,Ronco,Ronco2}).
A dendriform algebra is a family $(A,\prec,\succ)$ such that:
\begin{enumerate}
\item $A$ is a $K$-vector space and:
$$\begin{array}{rc|cl}
\prec: \left\{\begin{array}{rcl}
A\otimes A & \longrightarrow &A\\
a\otimes b& \longrightarrow & a\prec b,
\end{array}\right. &&&
\succ: \left\{\begin{array}{rcl}
A\otimes A & \longrightarrow &A\\
a\otimes b& \longrightarrow & a\succ b.
\end{array}\right. 
\end{array}$$
\item For all $a,b,c\in A$:
\begin{eqnarray}
\label{E1} (a\prec b)\prec c&=&a\prec(b\prec c+b\succ c),\\
\label{E2} (a\succ b)\prec c&=&a\succ(b\prec c),\\
\label{E3} (a\prec b+a\succ b)\succ c&=&a\succ(b\succ c).
\end{eqnarray}
\end{enumerate}}
\end{defi}

{\bf Remark.} If $A$ is a dendriform algebra, we put:
$$ m: \left\{\begin{array}{rcl}
A\otimes A & \longrightarrow &A\\
a\otimes b & \longrightarrow & ab=a\prec b+a\succ b.
\end{array}\right. $$
Then $(\ref{E1})+(\ref{E2})+(\ref{E3})$ is equivalent to the fact that $m$ is associative.
Hence, a dendriform algebra is a special (non unitary) associative algebra.\\

By duality, we obtain the notion of dendriform coalgebra:

\begin{defi}
\textnormal{
A dendriform coalgebra is a family $(C,\Delta_\prec,\Delta_\succ)$ such that:
\begin{enumerate}
\item $C$ is a $K$-vector space and:
$$\begin{array}{rc|cl}
\Delta_\prec: \left\{\begin{array}{rcl}
C & \longrightarrow &C\otimes C\\
a& \longrightarrow & \Delta_\prec(a)=a'_\prec \otimes a''_\prec,
\end{array}\right. &&&
\Delta_\succ: \left\{\begin{array}{rcl}
C & \longrightarrow &C\otimes C\\
a& \longrightarrow & \Delta_\succ(a)=a'_\succ \otimes a''_\succ.
\end{array}\right. 
\end{array}$$
\item for all $a\in C$:
\begin{eqnarray}
\label{E1'} (\Delta_\prec \otimes Id)\circ \Delta_\prec(a)&=&
(Id \otimes \Delta_\prec+Id \otimes \Delta_\succ)\circ \Delta_\prec(a),\\
\label{E2'} (\Delta_\succ \otimes Id)\circ \Delta_\prec(a)&=&
(Id \otimes \Delta_\prec)\circ \Delta_\succ(a),\\
\label{E3'} (\Delta_\prec \otimes Id+\Delta_\succ \otimes Id)\circ \Delta_\succ(a)&=&
(Id \otimes \Delta_\succ)\circ \Delta_\succ(a).
\end{eqnarray}
\end{enumerate}}
\end{defi}

{\bf Remarks.}
\begin{enumerate}
\item If  $C$ is a dendriform coalgebra, we put:
$$ \tdelta: \left\{\begin{array}{rcl}
C & \longrightarrow &C\otimes C\\
a & \longrightarrow & \tdelta(a)=\Delta_\prec(a)+\Delta_\succ(a)=a'\otimes a''.
\end{array}\right. $$
Then $(\ref{E1'})+(\ref{E2'})+(\ref{E3'})$ is equivalent
to the fact that $\tdelta$ is coassociative.
Hence, a dendriform coalgebra is a special (non counitary) coassociative coalgebra.
\item If $A=\bigoplus A_n$ is a $\mathbb{N}$-graded dendriform algebra, such 
that its homogeneous parts are finite-dimensional,
then $(A^{*g}, \prec^*, \succ^*)$ is a $\mathbb{N}$-graded 
dendriform coalgebra ($A^{*g}$ is the graded dual of $A$,
that is to say $A^{*g}=\bigoplus A_n^*\subseteq A^*$).
\item In the same way,  if $C$ is a $\mathbb{N}$-graded dendriform coalgebra, such 
that its homogeneous parts are finite-dimensional,
then $(C^{*g}, \Delta_\prec^*,\Delta_\succ^*)$ is a $\mathbb{N}$-graded dendriform algebra.
(In fact, for any dendriform coalgebra $C$,
the whole linear dual $C^*$ is a dendriform algebra).
\end{enumerate}

\begin{defi}
\textnormal{Let $A$ be a	dendriform coalgebra. We put:}
\begin{eqnarray*}
Prim(A)&=&\{a\in A\:/\:\tdelta(a)=0\},\\
Prim_\prec(A)&=&\{a\in A\:/\:\Delta_\prec(a)=0\},\\
Prim_\succ(A)&=&\{a\in A\:/\:\Delta_\succ(a)=0\},\\
Prim_{tot}(A)&=&Prim_\succ(A)\cap Prim_\prec(A)
\:=\:\{a\in A\:/\:\tdelta(a)=\Delta_\prec(a)=\Delta_\succ(a)=0\}.
\end{eqnarray*}
\end{defi}

\subsection{Dendriform and codendriform bialgebras}

\begin{defi}
\textnormal{(See \cite{Aguiar,Loday,Loday3,Ronco,Ronco2}).
A dendriform bialgebra is a family $(A,\prec,\succ,\tdelta)$ such that:
\begin{enumerate}
\item $(A,\prec,\succ)$ is a dendriform algebra.
\item $(A,\tdelta)$ is a coassociative (non counitary) coalgebra.
\item The following compatibilities are satisfied: for all $a,b \in A$,
\begin{eqnarray}
\label{E68} \tdelta(a\prec b)&=&a'b'\otimes a''\prec b''+a'\otimes a''\prec b
+a'b \otimes a''+b'\otimes a\prec b'' +b\otimes a,\\
\label{E57} \tdelta(a\succ b)&=&a'b'\otimes a''\succ b''+a'\otimes a''\succ b
+ab'\otimes b''+b'\otimes a\succ b''+a\otimes b.
\end{eqnarray}
\end{enumerate}}
\end{defi}

{\bf Remarks.}
\begin{enumerate}
\item  $(\ref{E68})+(\ref{E57})$ is equivalent to: for all $a \in A$,
\begin{equation}
\label{E5678} \tdelta(ab)=a'b'\otimes a''b''+a'\otimes a''b+ab'\otimes b''+a'b\otimes a''
+b'\otimes ab''+a\otimes b+b\otimes a.
\end{equation}
If $A$ is a dendriform bialgebra, we put $\overline{A}=A\oplus K$, which
is given a structure of associative algebra and coassociative coalgebra in the following way:
for all $a,b\in A$,
$$1.a=a,\hspace{1cm} a.1=a, \hspace{1cm} 1.1=1, \hspace{1cm} a.b= ab \mbox{ (product in $A$);}$$
$$\Delta(1)=1\otimes 1,\hspace{1cm} \Delta(a)=1\otimes a+a\otimes 1+\tdelta(a).$$
Then (\ref{E5678}) means that $\overline{A}$ is a bialgebra.
A dendriform bialgebra is then the augmentation
ideal of a special bialgebra.
\item Another interpretation of (\ref{E68}) and (\ref{E57}) will be given in section \ref{sect3}.\\
\end{enumerate}

By duality, we define the notion of codendriform bialgebra:
\begin{defi}
\textnormal{
A codendriform bialgebra is a family $(A,m,\Delta_\prec,\Delta_\succ)$ such that:
\begin{enumerate}
\item $(A,\Delta_\prec,\Delta_\succ)$ is a dendriform coalgebra.
\item $(A,m)$ is an associative (non unitary) algebra.
\item The following compatibilities are satisfied: for all $a,b \in A$,
\begin{eqnarray}
\label{E56} \Delta_\succ(ab)&=&a'b'_\succ\otimes a''b''_\succ+a'\otimes a''b
+ab'_\succ\otimes b''_\succ+b'_\succ\otimes ab''_\succ+a\otimes b,\\
\label{E78} \Delta_\prec(ab)&=& a'b'_\prec\otimes a''b''_\prec+a'b\otimes a''
+ab'_\prec\otimes b''_\prec+b'_\prec\otimes ab''_\prec+b\otimes a.
\end{eqnarray}
\end{enumerate}}
\end{defi}

{\bf Remarks.}
\begin{enumerate}
\item $(\ref{E56})+(\ref{E78})$ is equivalent to (\ref{E5678}).
Hence, if $A$ is a codendriform bialgebra, as before
$\overline{A}=A\oplus K$ is given a structure of bialgebra.
A codendriform bialgebra is then the augmentation
ideal of a special bialgebra.
\item If $A$ is a $\mathbb{N}$-graded codendriform bialgebra, 
such that its homogeneous parts are finite-dimensional,
then $(A^{*g}, \Delta_\prec^*, \Delta_\succ^*, m^*)$ is a $\mathbb{N}$-graded dendriform bialgebra.
\item In the same way,  if $A$ is a $\mathbb{N}$-graded dendriform bialgebra, such 
that its homogeneous parts are finite-dimensional,
then $(A^{*g}, \tdelta^*, \prec^*, \succ^*)$ is a $\mathbb{N}$-graded codendriform bialgebra.
\end{enumerate}

\subsection{Free dendriform algebras}

 Let us recall here the construction of  the Connes-Kreimer Hopf algebra of
planar decorated rooted trees (See \cite{Foissy2,Foissy3,Holtkamp} for more details).
It is a non-commutative version of the Connes-Kreimer Hopf algebra of rooted trees for Renormalisation
(\cite{Connes,Kreimer1,Kreimer2,Kreimer3}).

\begin{defi}
 \textnormal{
\begin{enumerate}
\item A {\it rooted tree} $t$ is a finite graph, without loops, with a special vertex called {\it root} of $t$. 
A {\it planar rooted tree} $t$ 
is a rooted tree with an imbedding in the plane.
The {\it weight} of $t$ is the number of its vertices. 
the set of planar rooted trees will be denoted by $\T$.
\item Let $\D$ be a nonempty set.
A planar rooted tree decorated by $\D$ is a planar tree
with an application from the set of its vertices into $\D$. 
The set of planar rooted trees decorated by $\D$ will be denoted by $\T^\D$.
\end{enumerate}}
\end{defi}

{\bf Examples.} \begin{enumerate}
\item Planar rooted trees with weight smaller than $5$:
$$\tun,\tdeux,\ttroisun,\ttroisdeux,\tquatreun, \tquatredeux,\tquatretrois,\tquatrequatre,\tquatrecinq,
\tcinqun, \tcinqdeux,\tcinqtrois,\tcinqquatre,\tcinqcinq,\tcinqsix,\tcinqsept,\tcinqhuit,
\tcinqneuf,\tcinqdix,\tcinqonze,\tcinqdouze,\tcinqtreize,\tcinqquatorze.$$
\item Planar rooted trees decorated by $\{a,b\}$ with weight smaller than $3$:
$$\tdun{$a$},\tdun{$b$}, \tddeux{$a$}{$a$},  \tddeux{$a$}{$b$}, 
\tddeux{$b$}{$a$}, \tddeux{$b$}{$b$},
\tdtroisun{$a$}{$a$}{$a$}, \tdtroisun{$b$}{$a$}{$a$}, 
\tdtroisun{$a$}{$b$}{$a$}, \tdtroisun{$a$}{$a$}{$b$}, 
\tdtroisun{$a$}{$b$}{$b$}, \tdtroisun{$b$}{$a$}{$b$}, 
\tdtroisun{$b$}{$b$}{$a$}, \tdtroisun{$b$}{$b$}{$a$},$$ 
$$\tdtroisdeux{$a$}{$a$}{$a$},\tdtroisdeux{$b$}{$a$}{$a$},
\tdtroisdeux{$a$}{$b$}{$a$},\tdtroisdeux{$a$}{$a$}{$b$},
\tdtroisdeux{$a$}{$b$}{$b$},\tdtroisdeux{$b$}{$a$}{$b$},
\tdtroisdeux{$b$}{$b$}{$a$},\tdtroisdeux{$b$}{$b$}{$b$}.$$
\end{enumerate}

The algebra $\h$ (denoted by $\h_{P,R}$ in \cite{Foissy2,Foissy3})
is the free associative (non commutative) $K$-algebra generated by the elements of $\T$.
Monomials in planar rooted trees in these algebra are called planar rooted forests.
The set of planar rooted forests will be denoted by $\F$.
Note that $\F$ is a basis of $\h$.

In the same way, if $\D$ is a nonempty set, the algebra $\h^\D$ 
(denoted by $\h_{P,R}^\D$ in \cite{Foissy2,Foissy3})
is the free associative (non commutative) $K$-algebra generated by the elements of $\T^\D$.
Monomials in planar rooted trees decorated by $\D$ in these algebra are called planar rooted forests
decorated by $\D$.
The set of planar rooted forests decorated by $\D$ will be denoted by $\F^\D$.
Note that $\F^\D$ is a basis of $\h^\D$.
if $\D$ is reduced to a single element, then  $\T^\D$ can be identified with $\T$
and $\h^\D$ can be identified with $\h$.\\

{\bf Examples.} 
\begin{enumerate}
\item Planar rooted forests of weight smaller than $4$:
$$1,\tun,\tun\tun,\tdeux,\tun\tun\tun,\tdeux\tun,\tun\tdeux,\ttroisun,\ttroisdeux,
\tun\tun\tun\tun,\tdeux\tun\tun,\tun\tdeux\tun,
\tun\tun\tdeux,\ttroisun\tun,\tun\ttroisun,\ttroisdeux\tun,\tun\ttroisdeux,
\tdeux\tdeux,\tquatreun,\tquatredeux,\tquatretrois,\tquatrequatre,\tquatrecinq.$$
\item Planar rooted forests decorated by $\{a,b\}$ of weight smaller than $2$:
$$1,\tdun{$a$},\tdun{$b$}, \tdun{$a$}\tdun{$a$}, \tdun{$b$}\tdun{$a$}, 
\tdun{$a$}\tdun{$b$}, \tdun{$b$}\tdun{$b$}, \tddeux{$a$}{$a$}, \tddeux{$b$}{$a$},
\tddeux{$a$}{$b$}, \tddeux{$b$}{$b$}.$$
\end{enumerate}

We now describe the Hopf algebra structure of $\h^\D$.
Let $t \in \T^\D$. 
An {\it admissible cut} of $t$ is a nonempty cut such that
every path in the tree meets at most one edge which is cut by $c$.
The set of admissible cut of $t$ is denoted by ${\cal A}dm(t)$. 
An admissible cut $c$ of $t$ sends $t$ to a couple $(P^c(t),R^c(t)) \in \F^\D \times \T^\D$, 
such that $R^c(t)$ is the connected component of the root of $t$
after the application of $c$, and $P^c(t)$ is the planar forest of the other
connected components (in the same order).
Moreover, if $c_v$ is the empty cut of $t$, we put $P^{c_v}(t)=1$ et $R^{c_v}(t)=t$. 
We define the {\it total cut} of $t$ as a cut $c_t$ such that  $P^{c_t}(t)=t$ and $R^{c_t}(t)=1$.
We then put ${\cal A}dm_*(t)={\cal A}dm(t)\cup\{c_v,c_t\}$. 

We now take $F \in \F^\D$, $F \neq 1$. 
There exists  $t_1, \ldots, t_n \in \T^\D$, such that $F=t_1\ldots t_n$.
An  {\it admissible cut} of $F$ is a $n$-uple $(c_1,\ldots,c_n)$ 
such that $c_i \in {\cal A}dm_*(t_i)$ for all $i$.
If all $c_i$'s are empty (resp. total), then $c$ is called the empty cut of $F$ 
(resp. the total cut of $F$). 
The set of admissible cuts of $F$
except the empty and the total cut is denoted by ${\cal A}dm(F)$.
The set of all admissible cuts of $F$ is denoted by ${\cal A}dm_*(F)$.
For $c=(c_1,\ldots,c_n) \in {\cal A}dm_*(F)$,
we put $P^c(F)= P^{c_1}(t_1) \ldots P^{c_n}(t_n)$ and  $R^c(F)= R^{c_1}(t_1) \ldots R^{c_n}(t_n)$.

The coproduct $\Delta: \h^\D \longrightarrow \h^\D \otimes \h^\D$ is defined in the following way:
for all $F \in \F^\D$,
$$\Delta(F)=F\otimes 1+1\otimes F+\sum_{c \in {\cal A}dm(F)} P^c(F)\otimes R^c(F)
=\sum_{c\in {\cal A}dm_*(F)} P^c(F)\otimes R^c(F).$$

{\bf Example.} If $a,b,c,d,e \in \D$:
\begin{eqnarray*}
\Delta\left( \tdun{$a$}\tdquatrequatre{$b$}{$c$}{$d$}{$e$}\right)&=&
 \tdun{$a$} \tdquatrequatre{$b$}{$c$}{$d$}{$e$}\otimes 1
+1\otimes \tdun{$a$} \tdquatrequatre{$b$}{$c$}{$d$}{$e$}
+\tdun{$a$}\otimes  \tdquatrequatre{$b$}{$c$}{$d$}{$e$}
 +\tdquatrequatre{$b$}{$c$}{$d$}{$e$}\otimes \tdun{$a$}
+\tdun{$a$}\tdun{$e$}\tdun{$d$} \otimes \tddeux{$b$}{$c$}
+\tdun{$e$}\tdun{$d$} \otimes \tdun{$a$}\tddeux{$b$}{$c$}\\
&&+\tdun{$a$}\tdtroisun{$c$}{$d$}{$e$}\otimes \tdun{$b$}
+\tdtroisun{$c$}{$d$}{$e$}\otimes \tdun{$a$}\tdun{$b$}
+\tdun{$a$}\tdun{$e$} \otimes \tdtroisdeux{$b$}{$c$}{$d$}
+\tdun{$e$} \otimes\tdun{$a$} \tdtroisdeux{$b$}{$c$}{$d$}
+\tdun{$a$}\tdun{$d$} \otimes \tdtroisdeux{$b$}{$c$}{$e$}
+\tdun{$d$} \otimes \tdun{$a$}\tdtroisdeux{$b$}{$c$}{$e$}.
\end{eqnarray*}

The counit $\varepsilon$ is given by:
$$\varepsilon: \left\{\begin{array}{rcl}
\h^\D& \longrightarrow &K\\
F\in \F^\D& \longrightarrow &\delta_{F,1}.
\end{array}\right. $$

We proved in \cite{Foissy3} that $\h^\D$ is isomorphic to
the free dendriform algebra generated by $\D$, which is described in \cite{Loday,Ronco}
in terms of planar binary trees. So the augmentation ideal $\A^\D$ of $\h^\D$ 
inherits a structure of dendriform algebra, also described in \cite{Foissy3}
with the help of another basis of $\h^\D$, introduced by duality.
Hence, $\A^\D$ is freely generated by the $\tdun{$d$}$'s, $d\in \D$,
as a dendriform algebra. 
Here is  an example of a computation of a product $\prec$ in $\A^\D$.
For all $x\in \A^\D$, $$\tdun{$d$}\prec x=B^+_d(x),$$ where
$B_d^+:\h^\D\longrightarrow \h^\D$ is the linear application
which send  a forest $t_1\ldots t_n$ to the planar decorated tree obtained by
grafting $t_1,\ldots,t_n$ on a common root decorated by $d$.
(This comes from the description of $\prec$ in terms of graftings in \cite{Foissy3}
and proposition 36 of \cite{Foissy2}).

As $\h^\D$ is self-dual (\cite{Foissy2}), $\A^\D$ is given a structure of codendriform bialgebra.
The description in \cite{Foissy3} of the left and right products in the dual basis of forests
allows us to describe this structure with the following definition:

\begin{defi}
\textnormal{
Let $F=t_1\ldots t_n \in \F^\D$, $F\neq 1$.
The set  $\A dm'_\prec(F)$ is the set of cuts $(c_1,\ldots,c_n)\in \A dm(F)$
such that $c_n$ is the total cut of $t_n$ if $F$ is not a single tree, and $\emptyset$ otherwise.
The set $\A dm'_\succ(F)$ is $\A dm(F)-\A dm'_\prec(F)$.}
\end{defi}

The dendriform coalgebra structure of $\A^\D$ is then given in the following way:
For all $F\in \F^\D-\{1\}$,
$$\Delta_\prec'(F)=\sum_{c\in \A dm_\prec(F)} P^c(F)\otimes R^c(F),\hspace{1cm}
\Delta_\succ'(F)=\sum_{c\in \A dm_\succ(F)} P^c(F)\otimes R^c(F).$$
The  product of $\A^\D$ is the product induced by the product of $\h^\D$.\\

{\bf Examples.}
\begin{enumerate}
\item If $t \in \T^\D$, $\Delta'_\prec(t)=0$.
\item If $a,b,c,d,e \in \D$.
\begin{eqnarray*}
\Delta'_\prec\left( \tdun{$a$}\tdquatrequatre{$b$}{$c$}{$d$}{$e$}\right)&=&
 \tdquatrequatre{$b$}{$c$}{$d$}{$e$}\otimes \tdun{$a$},\\
\Delta'_\succ\left( \tdun{$a$}\tdquatrequatre{$b$}{$c$}{$d$}{$e$}\right)&=&
\tdun{$a$}\otimes  \tdquatrequatre{$b$}{$c$}{$d$}{$e$}
+\tdun{$a$}\tdun{$e$}\tdun{$d$} \otimes \tddeux{$b$}{$c$}
+\tdun{$e$}\tdun{$d$} \otimes \tdun{$a$}\tddeux{$b$}{$c$}
+\tdun{$a$}\tdtroisun{$c$}{$d$}{$e$}\otimes \tdun{$b$}
+\tdtroisun{$c$}{$d$}{$e$}\otimes \tdun{$a$}\tdun{$b$}\\
&&+\tdun{$a$}\tdun{$e$} \otimes \tdtroisdeux{$b$}{$c$}{$d$}
+\tdun{$e$} \otimes\tdun{$a$} \tdtroisdeux{$b$}{$c$}{$d$}
+\tdun{$a$}\tdun{$d$} \otimes \tdtroisdeux{$b$}{$c$}{$e$}
+\tdun{$d$} \otimes \tdun{$a$}\tdtroisdeux{$b$}{$c$}{$e$}.
\end{eqnarray*}
\end{enumerate}

Moreover, $\h^\D$ can be graded.
A set $\D$ is said to be graded
when it is given an application $|.|:\D\longrightarrow \mathbb{N}$.
We denote $\D_n=\{d\in \D\:/\: |d|=n\}$.
We then put, for all $F\in \F^\D$:
$$|F|=\sum_{s\in vert(F)} |\mbox{decoration of $s$}|,$$ 
where $vert(F)$ is the set of vertices of $F$.
Then  $\h^\D$ is given a graded Hopf algebra structure
and $\A^\D$ is given a graded (co)dendriform bialgebra  structure by putting, for all $F\in \F^\D$,
$F$ homogeneous of degree $|F|$.
When $\D_0$ is empty and $\D_n$ is finite for all $n$, 
then  the homogeneous parts of $\A^\D$ are finite-dimensional and $(\A^\D)_0=(0)$.
Moreover, if those conditions occur, we have the following result:

\begin{prop}
\label{prop18}
We consider the following formal series:
$$D(X)=\sum_{n=1}^{+\infty}card(\D_n)X^n,\hspace{1cm}
R(X)=\sum_{n=0}^{+\infty}dim(\h^\D_n)X^n.$$
Then
$\displaystyle R(X)=\frac{1-\sqrt{1-4D(X)}}{2D(X)}$.
\end{prop}

{\bf Proof.} Similar to the proof of theorem 75 of \cite{Foissy4}. $\Box$

\section{Bidendiform bialgebras}

\label{sect2}

\subsection{Definition}

We now introduce the notion of bidendriform bialgebra. A bidendriform bialgebra
is both a dendriform bialgebra and a codendriform bialgebra, with some compatibilities.
\begin{defi}
\textnormal{
A bidendriform bialgebra is a family $(A,\prec,\succ,\Delta_\prec,\Delta_\succ)$ such that:
\begin{enumerate}
\item $(A,\prec,\succ)$ is a dendriform algebra.
\item $(A,\Delta_\prec,\Delta_\succ)$ is a dendriform  coalgebra.
\item The following compatibilities are satisfied: for all $a,b\in A$,
\begin{eqnarray}
\label{E5} \Delta_\succ (a\succ b)&=&a'b'_\succ\otimes a''\succ b''_\succ
+a'\otimes a''\succ b+b'_\succ\otimes a\succ b''_\succ+ab'_\succ \otimes b''_\succ+a\otimes b,\\
\label{E6} \Delta_\succ (a\prec b)&=&a'b'_\succ\otimes a''\prec b''_\succ
+a'\otimes a''\prec b+b'_\succ\otimes a\prec b''_\succ,\\
\label{E7} \Delta_\prec (a\succ b)&=& a'b'_\prec\otimes a''\succ b''_\prec
+ab'_\prec \otimes b''_\prec+b'_\prec\otimes a\succ b''_\prec,\\
\label{E8} \Delta_\prec (a\prec b)&=&a'b'_\prec\otimes a''\prec b''_\prec
+a'b \otimes a''+b'_\prec\otimes a\prec b''_\prec+b\otimes a.
\end{eqnarray}
\end{enumerate}}
\end{defi}

{\bf Remarks.}
\begin{enumerate}
\item  $(\ref{E6})+(\ref{E8})$ and $(\ref{E5})+(\ref{E7})$ are equivalent to 
(\ref{E68}) and (\ref{E57}), so a bidendriform bialgebra is a special dendriform bialgebra.
\item $(\ref{E5})+(\ref{E6})$ and $(\ref{E7})+(\ref{E8})$ are equivalent to (\ref{E56})
and (\ref{E78}), so a bidendriform bialgebra is a special codendriform bialgebra.
\item If $A$ is a graded bidendriform bialgebra, such 
that its homogeneous parts are finite-dimensional,
then $(A^{*g}, \Delta_\prec^*, \Delta_\succ^*, \prec^*,\succ^*)$
is also a graded bidendriform bialgebra, as the transposes of (\ref{E5}) and (\ref{E8}) are themselves, 
and (\ref{E6}) and (\ref{E7}) are transposes from each other.
\end{enumerate}

\subsection{Totally primitive elements}

We here give several results which will be useful to prove theorem \ref{theo17}. \\

{\it First part.} Let $A$ be a dendriform coalgebra. We define inductively:
$$\Delta_\prec^0=Id,\hspace{1cm}
\Delta_\prec^1=\Delta_\prec,\hspace{1cm}
\Delta_\prec^n=
(\Delta_\prec\otimes \underbrace{Id\otimes \ldots \otimes Id}_{n-1})\circ \Delta_\prec^{n-1}.$$
For all $n \in \mathbb{N}$, $\Delta_\prec^n:A\longrightarrow A^{\otimes(n+1)}$.

\begin{lemme}
\label{Lemme7}
For all $n \in \mathbb{N}^*$, $i\in \{2,\ldots,n\}$,
\begin{eqnarray*}
(\Delta_\prec \otimes \underbrace{Id\otimes \ldots \otimes Id}_{n-1})\circ \Delta_\prec^{n-1}
&=&\Delta_\prec^n,\\
(\underbrace{Id\otimes \ldots \otimes Id}_{i-1}\otimes \tdelta \otimes 
\underbrace{Id\otimes \ldots \otimes Id}_{n-i})\circ \Delta_\prec^{n-1}&=&\Delta_\prec^n.
\end{eqnarray*}
\end{lemme}

{\bf Proof.} The first equality is a consequence of the definition of $\Delta_\prec$.
Let us show the second one. By (\ref{E1'}):
\begin{eqnarray*}
&&(\underbrace{Id\otimes \ldots \otimes Id}_{i-1}\otimes \tdelta \otimes 
\underbrace{Id\otimes \ldots \otimes Id}_{n-i})\circ \Delta_\prec^{n-1}\\
&=&(\Delta_\prec\otimes \underbrace{Id\otimes \ldots \otimes Id}_{n-1})\circ
\ldots \circ (\Delta_\prec\otimes \underbrace{Id\otimes \ldots \otimes Id}_{n-i+2})\\
&&\circ (((Id \otimes \tdelta)\circ \Delta_\prec )\otimes 
\underbrace{Id\otimes \ldots \otimes Id}_{n-i})\circ
\Delta_\prec^{n-i}\\
&=&(\Delta_\prec\otimes \underbrace{Id\otimes \ldots \otimes Id}_{n-1})\circ
\ldots \circ (\Delta_\prec\otimes \underbrace{Id\otimes \ldots \otimes Id}_{n-i+2})\\
&&\circ (((\Delta_\prec \otimes Id)\circ \Delta_\prec )\otimes 
\underbrace{Id\otimes \ldots \otimes Id}_{n-i})\circ
\Delta_\prec^{n-i}\\
&=&(\Delta_\prec\otimes \underbrace{Id\otimes \ldots \otimes Id}_{n-1})\circ
\ldots \circ (\Delta_\prec\otimes \underbrace{Id\otimes \ldots \otimes Id}_{n-i+2})\\
&&\circ (\Delta_\prec \otimes \underbrace{Id\otimes \ldots \otimes Id}_{n-i+1})
 \circ (\Delta_\prec \otimes \underbrace{Id\otimes \ldots \otimes Id}_{n-i})\circ
\Delta_\prec^{n-i}\\
&=&\Delta_\prec^n.\:\Box
\end{eqnarray*}

\begin{lemme}
\label{lemme8}
Let $a\in A$, such that $\Delta_\prec^n(a)=0$. 
Then $\Delta_\prec^{n-1}(a)\in Prim_\prec(A)\otimes Prim(A)^{\otimes (n-1)}$.
\end{lemme}

{\bf Proof.} By lemma \ref{Lemme7},
$\Delta_\prec^{n-1}(a)$ vanishes under $\Delta_\prec \otimes Id^{\otimes(n-1)}$,
so belongs to $Prim_\prec(A)\otimes A^{\otimes (n-1)}$.
Moreover, if $i\geq 2$, $\Delta_\prec^{n-1}(a)$ vanishes under
$Id^{\otimes(i-1)}\otimes \tdelta \otimes Id^{\otimes(n-i)}$,
so belongs to $A^{\otimes(i-1)}\otimes Prim(A)\otimes A^{\otimes (n-i)}$. $\Box$\\

Suppose now that $A$ is a bidendriform bialgebra.
Let $a_1,\ldots,a_n \in A$. We define inductively:
\begin{eqnarray*}
\omega(a_1)&=&a_1,\\
\omega(a_1,a_2)&=&a_2\prec a_1,\\
\omega(a_1,\ldots,a_n)&=&a_n \prec \omega (a_1,\ldots,a_{n-1}).
\end{eqnarray*}

\begin{lemme}
\label{lemme9}
Let $a_1\in Prim_\prec(A)$ and $a_2,\ldots,a_n \in Prim(A)$. Let $k\in \mathbb{N}$. Then:
$$\Delta_\prec^k(\omega(a_1,\ldots,a_n))=
\sum_{1\leq i_1<i_2<\ldots<i_k<n} \omega(a_1,\ldots,a_{i_1})
\otimes \ldots \otimes \omega(a_{i_k+1},\ldots,a_n).$$
In particular:
$$  \left\{\begin{array}{rcl}
\Delta_\prec^{n-1}(\omega(a_1,\ldots,a_n))&=&a_1\otimes \ldots \otimes a_n,\\
\Delta_\prec^{k}(\omega(a_1,\ldots,a_n))&=&0\mbox{ if }k\geq n.
\end{array}\right. $$
\end{lemme}

{\bf Proof.} 
By induction on $k$. It is immediate for $k=0$.
Let us show the result for $k=1$ by induction on $n$. It is obvious for $n=1$.
We suppose that $n\geq 2$. As $a_n \in Prim(A)$, 
we have, by (\ref{E8}):
\begin{eqnarray*}
\Delta_\prec(\omega(a_1,\ldots,a_n))&=& \omega(a_1,\ldots,a_{n-1})'_\prec \otimes 
a_n \prec \omega(a_1,\ldots,a_{n-1})_\prec''
+\omega(a_1,\ldots,a_{n-1})\otimes a_n\\
&=&\sum_{1\leq i <n-2} \omega(a_1,\ldots,a_i)\otimes a_n \prec \omega(a_{i+1},\ldots,a_{n-1})
+\omega(a_1,\ldots,a_{n-1})\otimes a_n\\
&=&\sum_{1\leq i <n-2} \omega(a_1,\ldots,a_i)\otimes \omega(a_{i+1},\ldots,a_{n-1},a_n)
+\omega(a_1,\ldots,a_{n-1})\otimes a_n\\
&=&\sum_{1\leq i <n-1} \omega(a_1,\ldots,a_i)\otimes \omega(a_{i+1},\ldots,a_n).
\end{eqnarray*}

We suppose that the result is true at rank $k$. Then:
\begin{eqnarray*}
&&\Delta_\prec^{k+1}(\omega(a_1,\ldots,a_n))\\
&=&
\sum_{1\leq i_1<i_2<\ldots<i_k<n} \Delta_\prec(\omega(a_1,\ldots,a_{i_1}))
\otimes \omega(a_{i_1+1},\ldots a_{i_2})
\otimes \ldots \otimes \omega(a_{i_k+1},\ldots,a_n)\\
&=&\sum_{1\leq i_1<i_2<\ldots<i_{k+1}<n} \omega(a_1,\ldots,a_{i_1})
\otimes \ldots \otimes \omega(a_{i_{k+1}+1},\ldots,a_n). \:\Box\\
\end{eqnarray*}

For all $a \in A$, we put 
$N_\prec(a)=\inf \{n\in \mathbb{N}\:/\: \Delta_\prec^{n-1}(a)=0\}\in \mathbb{N}\cup\{+\infty\}$.

\begin{lemme}
\label{lemme10}
Let $a\in A$, such that $N_\prec(a)$ is finite. Then $a$ can be written as
a linear span of terms $\omega(a_1,\ldots,a_n)$, $n\leq N_\prec(a)$,
$a_1\in Prim_\prec(A)$, $a_2,\ldots,a_n \in Prim(A)$.
\end{lemme}

{\bf Proof.} By induction on $n=N_\prec(a)$. If $n=0$, then $a=0$.
If $n=1$, then $a\in Prim_\prec(A)$ and the result is obvious. Suppose $n\geq 2$.
By lemma \ref{lemme8}, we then put:
$$\Delta_\prec^{n-1}(a)=\sum a_1\otimes \ldots \otimes a_n 
\in Prim_\prec(A)\otimes Prim(A)^{\otimes (n-1)}.$$
By lemma \ref{lemme9}, $\Delta_\prec^{n-1}\left(a-\sum\omega(a_1,\ldots,a_n)\right)=0$.
The induction hypothesis applied to $a-\sum \omega(a_1,\ldots,a_n)$ gives the result. $\Box$ \\

{\it Second part.} Suppose that $A$ is a dendriform coalgebra.
We define inductively:
$$\tdelta^0=Id,\hspace{1cm}
\tdelta^1=\tdelta,\hspace{1cm}
\tdelta^n=
(\underbrace{Id\otimes \ldots \otimes Id}_{n-1}\otimes \tdelta)\circ \tdelta^{n-1}.$$
For all $n \in \mathbb{N}$, $\tdelta^n:A\longrightarrow A^{\otimes(n+1)}$.
Note that $\tdelta^n$ could be defined for any coassociative coalgebra $(A,\tdelta)$.

\begin{lemme}
\label{lemme11}
Let $a\in Prim_\prec(A)$. Then for all $n \in \mathbb{N}$,
$\tdelta^n(a)\in A^{\otimes n}\otimes Prim_\prec(A).$
\end{lemme}

{\bf Proof.} It is obvious if $n=0$. Suppose $n\geq 1$. By (\ref{E2'}):
\begin{eqnarray*}
(Id \otimes \Delta_\prec)\circ \tdelta&=& (Id \otimes \Delta_\prec)\circ \Delta_\prec
+(Id \otimes \Delta_\prec)\circ \Delta_\succ\\
&=&(Id \otimes \Delta_\prec)\circ \Delta_\prec
+(\Delta_\succ\otimes Id)\circ \Delta_\prec\\
&=&(Id \otimes \Delta_\prec+\Delta_\succ\otimes Id)\circ \Delta_\prec.
\end{eqnarray*}

By coassociativity of $\tdelta$, we have:
\begin{eqnarray*}
(Id^{\otimes n}\otimes \Delta_\prec)\circ \tdelta^n(a)
&=&(\tdelta^{n-1}\otimes Id \otimes Id)
\circ (Id \otimes \Delta_\prec)\circ \tdelta(a)\\
&=&(\tdelta^{n-1}\otimes Id \otimes Id)
\circ(Id \otimes \Delta_\prec+\Delta_\succ\otimes Id)\circ \Delta_\prec(a)\\
&=&0.
\end{eqnarray*}
Hence, $\tdelta^n(a)\in A^{\otimes n}\otimes Prim_\prec(A)$. $\Box$\\

\begin{lemme}
\label{lemme13}
Let $a\in A$, such that $\tdelta^n(a)=0$. 
Then $\tdelta^{n-1}(a)\in Prim(A)^{\otimes n}$.
Moreover, if $a \in Prim_\prec(a)$, then
$\tdelta^{n-1}(a)\in Prim(A)^{\otimes (n-1)}\otimes Prim_{tot}(A)$.
\end{lemme}

{\bf Proof.} By coassociativity of $\tdelta$,
for all $i\in \{1,\ldots, n\}$:
$$(Id^{\otimes(i-1)}\otimes \tdelta \otimes Id^{\otimes(n-i)})\circ \tdelta^{n-1}(a)
=\tdelta^n(a)=0,$$
so $\tdelta^{n-1}(a)\in Prim(A)^{\otimes n}$. The second assertion
comes from lemma \ref{lemme11}. $\Box$\\

Suppose now that $A$ is a bidendriform bialgebra.
Let $a_1,\ldots,a_n \in A$. We define inductively:
\begin{eqnarray*}
\omega'(a_1)&=&a_1,\\
\omega'(a_1,a_2)&=&a_1\succ a_2,\\
\omega'(a_1,\ldots,a_n)&=&\omega' (a_1,\ldots,a_{n-1})\succ a_n.
\end{eqnarray*}

\begin{lemme}
\label{lemme12}
Let $a_1,\ldots,a_n \in Prim(A)$. Let $k\in \mathbb{N}$. Then:
$$\tdelta^k(\omega'(a_1,\ldots,a_n))=
\sum_{1\leq i_1<i_2<\ldots<i_k<n} \omega'(a_1,\ldots,a_{i_1})
\otimes \ldots \otimes \omega'(a_{i_k+1},\ldots,a_n).$$
In particular:
$$  \left\{\begin{array}{rcl}
\tdelta^{n-1}(\omega'(a_1,\ldots,a_n))&=&a_1\otimes \ldots \otimes a_n,\\
\tdelta^{k}(\omega'(a_1,\ldots,a_n))&=&0\mbox{ if }k\geq n.
\end{array}\right. $$
Moreover, if $a_n \in Prim_{tot}(A)$, then $\omega'(a_1,\ldots,a_n)\in Prim_\prec(A)$.
\end{lemme}

{\bf Proof.} 
By induction on  $k$. It is immediate if $k=0$.
Let us show the result for $k=1$ by induction on $n$. It is obvious for $n=1$.
Suppose $n\geq 2$. As $a_n \in Prim(A)$, by (\ref{E57}):
\begin{eqnarray*}
\tdelta(\omega'(a_1,\ldots,a_n))&=& \omega'(a_1,\ldots,a_{n-1})'\otimes 
\omega'(a_1,\ldots,a_{n-1})''\succ a_n
+\omega(a_1,\ldots,a_{n-1})\otimes a_n,\\
&=&\sum_{1\leq i <n-2} \omega'(a_1,\ldots,a_i)\otimes \omega'(a_{i+1},\ldots,a_{n-1})\succ a_n
+\omega'(a_1,\ldots,a_{n-1})\otimes a_n\\
&=&\sum_{1\leq i <n-2} \omega'(a_1,\ldots,a_i)\otimes \omega'(a_{i+1},\ldots,a_{n-1},a_n)
+\omega'(a_1,\ldots,a_{n-1})\otimes a_n\\
&=&\sum_{1\leq i <n-1} \omega'(a_1,\ldots,a_i)\otimes \omega'(a_{i+1},\ldots,a_n).
\end{eqnarray*}

Suppose that the result is true at rank $k$. Then:
\begin{eqnarray*}
&&\tdelta^{k+1}(\omega'(a_1,\ldots,a_n))\\
&=&
\sum_{1\leq i_1<i_2<\ldots<i_k<n} \tdelta(\omega'(a_1,\ldots,a_{i_1}))
\otimes \omega'(a_{i_1+1},\ldots a_{i_2})
\otimes \ldots \otimes \omega'(a_{i_k+1},\ldots,a_n)\\
&=&\sum_{1\leq i_1<i_2<\ldots<i_{k+1}<n} \omega'(a_1,\ldots,a_{i_1})
\otimes \ldots \otimes \omega'(a_{i_{k+1}+1},\ldots,a_n). 
\end{eqnarray*}

Suppose $a_n \in Prim_\prec(A)$. By (\ref{E7}),
by putting $x=\omega'(a_1,\ldots,a_{n-1})$,
$\Delta_\prec(\omega'(a_1,\ldots,a_n))=\Delta_\prec(x\succ a_n)=0$.
So $\omega'(a_1,\ldots,a_n) \in Prim_\prec(A)$. $\Box$\\

For all  $a \in A$, we put 
$N(a)=\inf \{n\in \mathbb{N}\:/\: \tdelta^{n-1}(a)= 0\}\in \mathbb{N}\cup\{+\infty\}$.

\begin{lemme}
\label{lemme14}
Let $a\in Prim_\prec(A)$, such that $N(a)$ is finite.
Then $a$ can be written as a linear span of terms $\omega'(a_1,\ldots,a_n)$, $n\leq N(a)$,
$a_1,\ldots,a_{n-1} \in Prim(A)$, $a_n\in Prim_{tot}(A)$.
\end{lemme}

{\bf Proof.} Induction on $n=N(a)$. 
If $n=0$, then $a=0$. If $n=1$, then $a\in Prim_{tot}(A)=Prim(A)\cap Prim_\prec(A)$.
Suppose that $n\geq 2$ and put, by lemma \ref{lemme13}:
$$\tdelta^{n-1}(a)=\sum a_1\otimes \ldots \otimes a_n \in
Prim(A)^{\otimes (n-1)}\otimes Prim_{tot}(A).$$ We put $a'=a-\sum \omega'(a_1,\ldots,a_n)$. 
By lemma \ref{lemme12}, $a'\in Prim_\prec(A)$ and $N(a')<n$. Hence, $a'$ satisfies the induction hypothesis,
so the result is true for $a$. $\Box$\\

\subsection{Connected bidendriform bialgebras}

We prove in this paragraph that if $A$ is a connected (definition \ref{deficonnexe})
bidendriform bialgebra, it is generated by its totally primitive elements (theorem \ref{theo17}).\\

Let $C$ be a dendriform  coalgebra. We define inductively:
\begin{eqnarray*}
\P_C(0)&=&\{Id_C\} \subseteq \L(C),\\
\P_C(1)&=&\{\Delta_\prec,\:\Delta_\succ\}\subseteq \L(C,C^{\otimes 2}),\\
\P_C(n)&=& \left\{\left(Id^{\otimes (i-1)}\otimes \Delta_\prec \otimes Id^{\otimes (n-i)}\right)
\circ P\:/\: P\in \P_C(n-1),\: 1\leq i\leq n \right\}\\
&&\cup \left\{\left(Id^{\otimes (i-1)}\otimes \Delta_\succ \otimes Id^{\otimes (n-i)}\right)
\circ P\:/\: P\in \P_C(n-1),\: 1\leq i\leq n \right\}
\subseteq \L(C,C^{\otimes (n+1)}).\\
\end{eqnarray*}

\begin{defi}
\label{deficonnexe}
\textnormal{
\begin{enumerate}
\item Let $C$ be a dendriform coalgebra. It is said  {\it connected}
if, for all $a\in C$, there exists $n_a \in \mathbb{N}$, such that
for all $P\in \P_C(n_a)$, $P(a)=0$.
\item Let $C$ be a connected dendriform coalgebra. For all $a\in C$, we put:
$$deg_p(a)=\inf\{n\in \mathbb{N}\:/\: \forall P\in \P_C(n),\: P(a)=0\} \in \mathbb{N}.$$
\item Let $C$ be a connected dendriform coalgebra.
For all $n \in \mathbb{N}$, we put:
$$C^{\leq n}=\{a\in C\:/\: deg_p(a)\leq n\}.$$
Then $C^{\leq 0}=(0)$, $C^{\leq 1}=Prim_{tot}(C)$,
and $(C^{\leq n})_{n\in \mathbb{N}}$ is a increasing filtration.
\end{enumerate}}
\end{defi}

{\bf Remarks.}
\begin{enumerate}
\item For all $n \in \mathbb{N}$, $\Delta_\prec^n\in \P_C(n)$ et 
$\tdelta^n \in Vect(\P_C(n))$. Hence, if $C$ is connected,
then for all $a\in C$, $N_\prec(a)$ and $N(a)$ are finite and smaller than $deg_p(a)$.
\item Let $C$ be a $\mathbb{N}$-graded dendriform coalgebra, such that $C_0=(0)$ 
(the homogeneous parts of $C$ may not be finite-dimensional).
 Then $C$ is connected, as, for all $a\in C_n$, for all $P\in \P(n)$: 
$$P(a)\in \bigoplus_{k_1+\ldots+k_{n+1}=n} C_{k_1}\otimes \ldots \otimes C_{k_{n+1}}=(0).$$
Moreover, for all $a \in C$, $deg_p(a)\leq |a|$, where $|a|$ is the degree of $a$
for the gradation of $C$.
\end{enumerate}

\begin{lemme}
\label{le15}
Let $C$ be a connected dendriform coalgebra.
For all $n\geq 1$:
\begin{eqnarray*}
\Delta_\prec(C^{\leq n})&\subseteq & C^{\leq n-1}\otimes C^{\leq n-1},\\
\Delta_\succ(C^{\leq n})&\subseteq & C^{\leq n-1}\otimes C^{\leq n-1},\\
\tdelta(C^{\leq n})&\subseteq & C^{\leq n-1}\otimes C^{\leq n-1}.
\end{eqnarray*}
\end{lemme}

{\bf Proof.} Let $a\in C^{\leq n}$. Let $P\in \P_C(n-1)$. Then
$(P\otimes Id)\circ \Delta_\prec \in \P_C(n)$, so: 
$$(P\otimes Id)\circ \Delta_\prec(a)=0.$$
Hence, $\Delta_\prec(a)\in C^{\leq n-1}\otimes C$.
In the same way, we obtain $\Delta_\prec(a)\in C\otimes C^{\leq n-1}$
by considering $(Id \otimes P)\circ \Delta_\prec$,
so $\Delta_\prec(a)\in C^{\leq n-1}\otimes C^{\leq n-1}$.
The procedure is the same for $\Delta_\succ(a)$,
and the result for $\tdelta(a)$ is obtained by addition. $\Box$\\

The following lemma is now immediate:
\begin{lemme}
\label{lemmeutile}
Let $C$ be a connected dendriform coalgebra and $k,n\in \mathbb{N}^*$. Then:
$$\Delta_\prec^k\left(C^{\leq n}\right)\subseteq \left( C^{\leq n-1}\right)^{\otimes (k+1)},\hspace{1cm}
\tdelta^k\left(C^{\leq n}\right)\subseteq \left( C^{\leq n-1}\right)^{\otimes (k+1)}.$$
\end{lemme}

\begin{theo}
\label{theo17}
Let $A$ be a connected (as a dendriform coalgebra)
bidendriform bialgebra. Then $A$ is generated (as a dendriform algebra) by $Prim_{tot}(A)$.
\end{theo}

{\bf Proof.}
Let $B$ be the dendriform subalgebra of $A$ generated by $Prim_{tot}(A)$.
Let $a\in A$. We denote $deg_p(a)=n$. Let us show that $a \in B$ by induction on $n$.
If $n=0$, then $a=0 \in B$. If $n=1$,
then $a\in Prim_{tot}(A)\subseteq B$. Suppose $n\geq 2$.
As $A$ is connected, by remark $1$ after definition \ref{deficonnexe}, 
$N_\prec(a)=k$ is finite and smaller than $n$.
By lemma \ref{lemme10}, we can suppose $a=\omega(a_1,\ldots,a_k)$,  $k\leq n$, $a_1\in Prim_\prec(A)$,
$a_2,\ldots,a_k \in Prim(A)$.
We have two different cases.
\begin{enumerate}
\item If $k\geq 2$, by lemma \ref{lemmeutile},
$\Delta_\prec^{k-1}(\omega(a_1,\ldots,a_k))=a_1\otimes \ldots \otimes a_k
 \in \left(A^{\leq n-1}\right)^{\otimes k}$,
so, for all $i$, $deg_p(a_i)<n$, so $a_i \in B$. Hence, $a\in B$.
\item If $k=1$, then $a\in Prim_\prec(A)$. 
As $A$ is connected, by remark $1$ after definition \ref{deficonnexe}, 
$N(a)=l$ is finite ans smaller than $n$.
By lemma \ref{lemme14}, we can suppose $a=\omega'(b_1,\ldots,b_l)$, 
$l\leq n$, $b_1,\ldots,b_{l-1} \in Prim(A)$, $b_l\in Prim_{tot}(A)$.
We have two different cases.
\begin{enumerate}
\item If $l\geq 2$, by lemma \ref{lemmeutile},
$\tdelta^{l-1}(\omega'(b_1,\ldots,b_l))=b_1\otimes \ldots \otimes b_l
 \in \left(A^{\leq n-1}\right)^{\otimes l}$,
so, for all $i$, $deg_p(b_i)<n$, so $b_i \in B$. Hence, $a\in B$.
\item If $l=1$, then $a=b_1\in Prim_{tot}(A)\subseteq B$. $\Box$
\end{enumerate}
\end{enumerate}

\subsection{Projections on $Prim_\prec(A)$ and $Prim_{tot}(A)$}

We here define an eulerian idempotent for connected bidendriform bialgebras.\\

Let $A$ be a bidendriform bialgebra.
We put:
$$A^{D2}=Vect(x\prec y,\:x\succ y\:/\:x,y\in A).$$

We define $m_\prec^n:A^{\otimes n} \longrightarrow A$
inductively:
\begin{eqnarray*}
m_\prec^1(a_1)&=&a_1,\\
m_\prec^2(a_1\otimes a_2)&=&a_2\prec a_1,\\
m_\prec^n(a_1\otimes \ldots \otimes a_n)&=&m_\prec^{n-1}(a_2\otimes \ldots\otimes a_n) \prec a_1.
\end{eqnarray*}

\begin{prop}
\label{proppro}
Let $A$ be a bidendriform bialgebra such that for all $a\in A$,
$N_\prec(a)$ is finite.
We consider the following application:
$$T_1: \left\{\begin{array}{rcl}
A& \longrightarrow &A\\
a& \longrightarrow &\displaystyle \sum_{k=1}^{+\infty} 
(-1)^{k+1}m_\prec^k \circ \Delta_\prec^{k-1}(a).
\end{array}\right. $$
Then $T_1$ is a projection on $Prim_\prec(A)$. Moreover, for all $a\in A$,
$T_1(a)=a+A^{D2}$.
\end{prop}

{\bf Proof.} Remark first that $T_1$ is well defined, 
as for all $a\in A$, $\Delta_\prec^{k-1}(a)=0$ for $k$ great enough.
Let us show that $T_1(a)\in Prim_\prec(A)$ for all $a$.
By lemma \ref{lemme10}, we can suppose $a=\omega(a_1,\ldots,a_n)$, avec $a_1\in Prim_\prec(A)$,
$a_2,\ldots,a_n \in Prim(A)$. If $n=1$, then $T_1(a)=a_1\in Prim_\prec(A)$.
Suppose $n\geq 2$.

We consider the following binary trees: for all $k \in \mathbb{N}$,
$$t_k^{(d)}=\overbrace{
\begin{picture}(15,30)(-7,0)
\put(-3.5,4){$\vee$}
\put(0,0){\line(0,1){5}}
\put(-.5,10){$\vee$}
\put(6,17.5){\tiny .}
\put(6.5,19){\tiny .}
\put(7,20.5){\tiny .}
\put(5,22){$\vee$}
\put(0.5,0){\tiny $\prec$}
\put(3,6){\tiny $\prec$}
\put(9,18){\tiny $\prec$}
\end{picture}}^{\mbox{\tiny $k$ leaves}}\:;\:
t_k^{(g)}=\overbrace{\begin{picture}(10,30)(-7,0)
\put(-3.5,4){$\vee$}
\put(0,0){\line(0,1){5}}
\put(-6,10){$\vee$}
\put(-6.5,17.5){\tiny .}
\put(-7,19){\tiny .}
\put(-7.5,20.5){\tiny .}
\put(-11,22){$\vee$}
\put(-6,0){\tiny $\prec$}
\put(-8.5,6){\tiny $\prec$}
\put(-14.5,18){\tiny $\prec$}
\end{picture}}^{\mbox{\tiny $k$ leaves}}.
$$

Then $m_\prec^k$ and $\omega$ can be  graphically represented in the following way:
$$  \left\{\begin{array}{rcl}
\omega(x_1,\ldots,x_k)&=&t_k^{(d)}.x_k\otimes \ldots \otimes x_1, \\
m_\prec^k(x_1\otimes \ldots \otimes x_k)&=&t_k^{(g)}.x_k\otimes \ldots \otimes x_1.
\end{array}\right. $$
By lemma \ref{lemme9}, we have:
$$\Delta^{k-1}_\prec(\omega(a_1,\ldots,a_n))=
\sum_{n_1+\ldots+n_k=n}t_{n_1}^{(d)}.(a_{n_1}\otimes \ldots \otimes a_1)\otimes \ldots
\otimes t_{n_k}^{(d)}.(a_n \otimes \ldots \otimes a_{n_1+\ldots +n_{k-1}+1}).
$$
(The $n_i$'s are positive, non zero integers). Then:
$$T_1(\omega(a_1,\ldots,a_n)) = \sum_{k=1}^{+\infty}\: \sum_{n_1+\ldots+n_k=n}
(-1)^{k+1}t_{n_1,\ldots,n_k}.(a_n\otimes \ldots \otimes a_1),$$
where the  $n_i$'s are positive, non zero integers and
 $t_{n_1,\ldots,n_k}$ is the following tree:
$$\begin{picture}(18,50)(0,0)
\put(0,0){\line(0,1){5}}
\put(0,5){\line(-1,1){30}}
\put(0,5){\line(1,1){10}}
\put(10,16){$t_{n_k}^{(d)}$}
\put(-20,25){\line(1,1){10}}
\put(-10,38){$t_{n_2}^{(d)}$}
\put(-35,38){$t_{n_1}^{(d)}$}
\put(2,23){\tiny .}
\put(0,25){\tiny .}
\put(-2,27){\tiny .}
\put(-10,0){$\prec$}
\put(-30,20){$\prec$}
\end{picture}.$$
As  $t_{1,n_2,\ldots,n_k}=t_{1+n_2,\ldots,n_k}$, we obtain:
\begin{eqnarray*}
T_1(\omega(a_1,\ldots,a_n)) &=&\sum_{k=1}^{+\infty}\:
 \sum_{\stackrel{n_1+\ldots+n_k=n,}{n_1\geq 2}}\:
(-1)^{k+1}t_{n_1,\ldots,n_k}.(a_n\otimes \ldots \otimes a_1)\\
&&+\sum_{k=1}^{+\infty}\: \sum_{n_1+\ldots+n_k=n-1}\:
(-1)^{k+2}t_{1,n_1,\ldots,n_k}.(a_n\otimes \ldots \otimes a_1)\\
&=&\sum_{k=1}^{+\infty}\:
 \sum_{\stackrel{n_1+\ldots+n_k=n,}{n_1\geq 2}}\:
(-1)^{k+1}t_{n_1,\ldots,n_k}.(a_n\otimes \ldots \otimes a_1)\\
&&-\sum_{k=1}^{+\infty}\:
 \sum_{\stackrel{n'_1+\ldots+n_k=n,}{n'_1\geq 2}}\:
(-1)^{k+1}t_{n'_1,\ldots,n_k}.(a_n\otimes \ldots \otimes a_1)\:\mbox{ (where $n_1'=1+n_1$)}\\
&=&0.
\end{eqnarray*}
Hence, $T_1(a)\in Prim_\prec(A)$ for all $a\in A$.
Moreover, if $a\in Prim_\prec(A)$, $T_1(a)=a$.
So $T_1$ is a projection on $Prim_\prec(A)$. Finally, for all $a \in A$, we have:
$$T_1(a)=a+\sum_{k=2}^{+\infty} (-1)^{k+1}m_\prec^k \circ \Delta_\prec^{k-1}(a)=a+A^{D2}.\: \Box$$

We define $m_\succ^n:A^{\otimes n} \longrightarrow A$ inductively:
\begin{eqnarray*}
m_\succ^1(a_1)&=&a_1,\\
m_\succ^2(a_1\otimes a_2)&=&a_1\succ a_2,\\
m_\succ^n(a_1\otimes \ldots \otimes a_n)&=&a_1 \succ 
m_\succ^{n-1}(a_2\otimes \ldots\otimes a_n).
\end{eqnarray*}
With lemmas \ref{lemme11} and \ref{lemme12}, we can show
the following proposition in the same way as proposition \ref{proppro}:

\begin{prop}
Let $A$ be a bidendriform bialgebra such that for all $a\in Prim_\prec(A)$,
$N(a)$ is finite.
We consider the following application:
$$T_2: \left\{\begin{array}{rcl}
Prim_\prec(A)& \longrightarrow &A\\
a& \longrightarrow &\displaystyle \sum_{k=1}^{+\infty} (-1)^{k+1}m_\succ^k \circ \tdelta^{k-1}(a).
\end{array}\right. $$
Then $T_2$ is a projector from $Prim_\prec(A)$ into $Prim_{tot}(A)$. Moreover, for all $a\in A$,
$T_2(a)=a+A^{D2}$.
\end{prop}

\begin{cor}
\label{corutoi}
Let $A$ be a bidendriform bialgebra such that for all $a\in A$,
$N_\prec(a)$ is finite and for all $A\in Prim_\prec(A)$, $N(a)$ is finite (for example, $A$ is connected).
We consider the application $T=T_2\circ T_1:A \longrightarrow A$.
Then $T$ is a projection from $A$ into $Prim_{tot}(A)$.
Moreover, for all $a\in A$, $T(a)=a+A^{D2}$.
So $A=Prim_{tot}(A)+A^{D2}$.
\end{cor}

{\bf Proof.} Immediate, by composition. $\Box$

\begin{cor}
\label{cor15}
Let $A$ be a $\mathbb{N}$-graded bidendriform bialgebra, such that the homogeneous part
of $A$ are finite-dimensional and $A_0=(0)$. 
Then $A=Prim_{tot}(A)\oplus A^{D2}$.
\end{cor}

{\bf Proof.} Then $A$ is connected (remark $2$ after definition \ref{deficonnexe}).
By corollary \ref{corutoi}, $A=Prim_{tot}(A)+A^{D2}$.
In the same way, the bidendriform bialgebra $A^{*g}$ is connected.
We then have $A^{*g}=Prim_{tot}(A^{*g})+(A^{*g})^{D2}$. By taking the orthogonal:
$$(0)=(A^{*g})^\perp=Prim_{tot}(A^{*g})^\perp \cap ((A^{*g})^{D2})^\perp=A^{D2}\cap Prim_{tot}(A).$$
So $A=Prim_{tot}(A)\oplus A^{D2}$. $\Box$

\section{Tensor product and dendriform modules}

\label{sect3}

\subsection{Tensor product of dendriform algebras}

We here show how the category of dendriform algebras can be given a structure of tensor category.
As dendriform algebras are not objects with unit, we have to extend the usual tensor product 
in order to obtain a copy of $A$ and $B$ in the tensor product of $A$ and $B$.

\begin{defi}
\textnormal{Let $A,B$ be two vector spaces. Then:}
$$A\totimes B=(A\otimes B) \oplus (K\otimes B)\oplus (A \otimes K).$$
\end{defi}

\begin{prop}
Let $A,B$ be two dendriform algebras. Then $A\totimes B$ is given a structure of dendriform
algebra in the following way: for all $a_1,a_2 \in A$, $b_1,b_2 \in B$,
$$\begin{array}{rclc|crcl}
(a_1\otimes b_1)\prec (a_2\otimes b_2)&=&a_1.a_2 \otimes b_1 \prec b_2,&&&
(a_1\otimes b_1)\succ (a_2\otimes b_2)&=&a_1.a_2 \otimes b_1 \succ b_2,\\
(a_1\otimes b_1)\prec (a_2\otimes 1)&=&a_1.a_2 \otimes b_1,&&&
(a_1\otimes b_1)\succ (a_2\otimes 1)&=&0,\\
(a_1\otimes b_1)\prec (1\otimes b_2)&=&a_1 \otimes b_1 \prec b_2,&&&
(a_1\otimes b_1)\succ (1\otimes b_2)&=&a_1 \otimes b_1 \succ b_2,\\
(a_1\otimes 1)\prec (a_2\otimes b_2)&=&0,&&&
(a_1\otimes 1)\succ (a_2\otimes b_2)&=&a_1.a_2 \otimes b_2,\\
(a_1\otimes 1)\prec (a_2\otimes 1)&=&a_1\prec a_2 \otimes 1,&&&
(a_1\otimes 1)\succ (a_2\otimes 1)&=&a_1\succ a_2 \otimes 1,\\
(a_1\otimes 1)\prec (1\otimes b_2)&=&0,&&&
(a_1\otimes 1)\succ (1\otimes b_2)&=&a_1 \otimes b_2,\\
(1 \otimes b_1)\prec (a_2 \otimes b_2)&=&a_2 \otimes b_1 \prec b_2,&&&
(1 \otimes b_1)\succ (a_2 \otimes b_2)&=&a_2 \otimes b_1 \succ b_2,\\
(1\otimes b_1)\prec (a_2\otimes 1)&=&a_2 \otimes b_1,&&&
(1\otimes b_1)\succ (a_2\otimes 1)&=&0,\\
(1\otimes b_1)\prec (1\otimes b_2)&=&1 \otimes b_1 \prec b_2;&&&
(1\otimes b_1)\succ (1\otimes b_2)&=&1 \otimes b_1 \succ b_2.
\end{array}$$
\end{prop}

{\bf Proof.} Direct computations. $\Box$\\

{\bf Remark.} $A\otimes K$ is a dendriform subalgebra of $A\totimes B$ which is
isomorphic to $A$ and $K\otimes B$ is a dendriform subalgebra of $A\totimes B$ which is
isomorphic to $B$. 

\begin{prop}
\begin{enumerate}
\item Let $A,B,C$ be unitary dendriform algebras. Then the following application
is an isomorphism of dendriform algebras:
$$  \left\{\begin{array}{rcl}
(A\totimes B)\totimes C & \longrightarrow &A \totimes (B\totimes C)\\
(a\otimes b)\otimes c& \longrightarrow & a\otimes (b\otimes c).
\end{array}\right. $$
\item Let $A,A',B,B'$ be unitary dendriform algebras and $\phi:A\longrightarrow A'$,
$\psi:B\longrightarrow B'$ be morphisms of unitary dendriform algebras.
We then define:
$$\phi \totimes \psi: \left\{\begin{array}{rcl}
A\totimes B& \longrightarrow &A'\totimes B'\\
a\otimes b&\longrightarrow & \phi(a)\otimes \psi(b),\\
a\otimes 1&\longrightarrow & \phi(a)\otimes 1,\\
1\otimes b&\longrightarrow & 1\otimes \psi(b),
\end{array}\right. $$
for all $a\in A$, $b\in B$.
Then $\phi \totimes \psi$ is a morphism of dendriform algebras. 
\end{enumerate}
In other terms, the category of unitary dendriform algebras is a tensor category
with $\totimes$.
\end{prop}

{\bf Proof.} Direct computations. $\Box$\\

We can now reformulate the axioms of bidendriform bialgebras.
Let $(A,\prec,\succ)$ be a dendriform algebra 
and let $\Delta_\prec,\Delta_\succ:A\longrightarrow A\otimes A$.
We put $\tdelta=\Delta_\prec+\Delta_\succ$. We then define:
$$\begin{array}{cl}
\overline{\Delta}_\prec:& \left\{\begin{array}{rcl}
A\oplus K& \longrightarrow &A\totimes A\\
a& \longrightarrow &\Delta_\prec (a)+a\otimes 1,\\
1&\longrightarrow & 0,\\
\end{array}\right. \\ \\
\overline{\Delta}_\succ:& \left\{\begin{array}{rcl}
A\oplus K& \longrightarrow &A\totimes A\\ 
a& \longrightarrow &\Delta_\succ (a)+1\otimes a,\\
1&\longrightarrow & 0,\\
\end{array}\right. \\ \\
\Delta:& \left\{\begin{array}{rcl}
A\oplus K& \longrightarrow &A\totimes A\\
a& \longrightarrow &\tdelta(a)+a\otimes 1+1\otimes a,\\
1&\longrightarrow & 0.
\end{array}\right. 
\end{array}$$
Note that $\Delta=\overline{\Delta}_\prec+\overline{\Delta}_\succ$.
We can easily prove the following assertions:
\begin{enumerate}
\item (\ref{E1'})-(\ref{E3'}) for $\Delta_\prec$ and $\Delta_\succ$
are equivalent to (\ref{E1'})-(\ref{E3'}) for $\overline{\Delta}_\prec$ 
and $\overline{\Delta}_\succ$.
\item (\ref{E5})-(\ref{E8}) is equivalent to: for all $a,b\in A$,
$$\begin{array}{rclccrcl}
\overline{\Delta}_\succ(a\succ b)&=&\Delta(a) \succ \overline{\Delta}_\succ(b),&&&
\overline{\Delta}_\prec(a\succ b)&=&\Delta(a) \succ \overline{\Delta}_\prec(b),\\
\overline{\Delta}_\succ(a\prec b)&=&\Delta(a) \prec \overline{\Delta}_\succ(b),&&&
\overline{\Delta}_\prec(a\prec b)&=&\Delta(a) \prec \overline{\Delta}_\prec(b).
\end{array}$$
\item (\ref{E57})-(\ref{E68}) is equivalent to: for all $a,b \in A$,
$$\begin{array}{rclccrcl}
\Delta(a\succ b)&=&\Delta(a) \succ \Delta(b),&&&
\Delta(a\prec b)&=&\Delta(a) \prec \Delta(b).
\end{array}$$
In other terms, (\ref{E57})-(\ref{E68}) is equivalent to  $\Delta:A\longrightarrow A\totimes A$
is a morphism of dendriform algebras.
\end{enumerate}

\subsection{Dendriform modules}

(See for example \cite{Markl} for complements on operads).
We denote by  $\P_{dend}=(\P_{dend}(n))_{n\in \mathbb{N}^*}$
 the operad of dendriform algebras (we consider here non-$\Sigma$-operads,
that is to say there is no action of the symmetric groups).
In other words, $\P_{dend}$ is the operad generated by $\prec$ and $\succ \in \P_{dend}(2)$
and the following relations:
\begin{eqnarray*}
\prec\circ (\prec, I)&=&\prec \circ(I,\prec+\succ),\\
\prec \circ ( \succ,I)&=&\succ \circ (I,\prec),\\
\succ \circ (\prec+\succ,I)&=&\succ\circ (I,\succ).
\end{eqnarray*}

Let $A$ be a dendriform algebra. A dendriform module over $A$ is a vector space $M$
together with applications, for all $n \in \mathbb{N}^*$:
$$ \left\{\begin{array}{rcl}
\P_{dend}(n)\otimes A^{\otimes (n-1)}\otimes M& \longrightarrow &M\\
p\otimes (a_1\otimes \ldots \otimes a_{n-1})\otimes m& \longrightarrow &
p.(a_1\otimes \ldots \otimes a_{n-1}\otimes m),
\end{array}\right. $$
satisfying the same associativity relations and unit relation as those for $\P_{dend}$-algebras.
In other terms, a dendriform module over $A$ is a vector space $M$
together with two applications:
$$\dashv: \left\{\begin{array}{rcl}
A\otimes M& \longrightarrow &M\\
a\otimes m& \longrightarrow &a\dashv m,
\end{array}\right. \hspace{1cm}
\vdash: \left\{\begin{array}{rcl}
A\otimes M& \longrightarrow &M\\
a\otimes m& \longrightarrow &a\vdash m,
\end{array}\right. $$
with the following compatibilities: for all $a,b\in A$, $m\in M$,
\begin{eqnarray}
\label{E16} (a\prec b)\dashv m&=&a\dashv(b\dashv m+b\vdash m),\\
\label{E17} (a\succ b)\dashv m&=&a\vdash(b\dashv m),\\
\label{E18} (a\prec b+a\succ b)\vdash m&=&a\vdash(b\vdash m).
\end{eqnarray}
(We have $a\dashv m=\prec.(a\otimes m)$ and $a\vdash m=\succ.(a\otimes m)$).\\

{\bf Remarks.}
\begin{enumerate}
\item  If $M$ is a dendriform module over $A$, then $\dashv+\vdash$
gives $M$ a structure of left module over the associative algebra $(A,\prec+\succ)$ 
by (\ref{E16})+(\ref{E17})+(\ref{E18}).
This action will be denoted by $a.m=a\dashv m+a\vdash m$.
\item $A$ is a dendriform module over itself with $\dashv=\prec$
and $\vdash=\succ$.
\item Suppose that $A$ is a dendriform bialgebra. Then, as $\Delta$ is a morphism of dendriform
algebras, $A\totimes A$ is given a structure of a dendriform module over $A$ with,
for all $a\in A$, $\sum b\otimes c \in A\totimes A$:
$$a\dashv\left(\sum b\otimes c\right)=\sum \Delta(a)\prec(b\otimes c),\hspace{1cm}
a\vdash\left(\sum b\otimes c\right)=\sum \Delta(a)\succ(b\otimes c).$$
Moreover, $\Delta$ is a morphism of dendriform modules.
\item Let $(A,\prec,\succ,\tdelta)$ be a dendriform bialgebra 
and let $\Delta_\prec,\Delta_\succ:A\longrightarrow A\otimes A$.
Observe that (\ref{E5})-(\ref{E8}) is equivalent to: 
$\overline{\Delta}_\prec,\overline{\Delta}_\succ$ are morphisms of dendriform modules.
\end{enumerate}

\begin{prop}
\label{propo29}
Let $A$ be a free dendriform algebra generated by a subspace $V$ and let $M$ be
a dendriform module over $A$. For any linear application $\phi:V\longrightarrow M$,
there exists a unique morphism of dendriform modules $\Phi:A\longrightarrow M$,
such that $\Phi_{\mid V}=\phi$.
\end{prop}

{\bf Proof.}

{\it Unicity.} Because $V$ generates $A$ as a dendriform algebra,
and hence as a dendriform module.

{\it Existence.} As $A$ is freely generated by $V$,
we can suppose that $\displaystyle A=\bigoplus_{n\geq 1} \P_{dend}(n) \otimes V^{\otimes n}$.
We then define:
$$\Phi: \left\{\begin{array}{rcl}
A& \longrightarrow &M\\
p\otimes (v_1\otimes \ldots \otimes v_n)& \longrightarrow &
p.(v_1\otimes \ldots \otimes v_{n-1} \otimes \phi(v_n))
\end{array}\right. $$
for all $p\in \P_{dend}(n)$, $v_1,\ldots,v_n \in V$.
Then $\Phi$ fits the asked conditions. $\Box$\\

\begin{prop}
\label{propo30}
Let $(A,\prec,\succ,\tdelta)$ be a dendriform bialgebra, with applications
$\Delta_\prec,\Delta_\succ:A\longrightarrow A\otimes A$. We put:
$$\begin{array}{cl}
\overline{\Delta}_\prec:& \left\{\begin{array}{rcl}
A& \longrightarrow &A\totimes A\\
a& \longrightarrow &\Delta_\prec (a)+a\otimes 1,\\
\end{array}\right. \\ \\
\overline{\Delta}_\succ:& \left\{\begin{array}{rcl}
A& \longrightarrow &A\totimes A\\ 
a& \longrightarrow &\Delta_\succ (a)+1\otimes a,\\
\end{array}\right. \\ \\
\Delta:& \left\{\begin{array}{rcl}
A& \longrightarrow &A\totimes A\\
a& \longrightarrow &\tdelta(a)+a\otimes 1+1\otimes a.
\end{array}\right. 
\end{array}$$
We suppose the following conditions:
\begin{enumerate}
\item $(A,\Delta_\prec,\Delta_\succ)$ satisfies
relations (\ref{E1'}), (\ref{E2'}) and (\ref{E3'}) on a set
of generators of the dendriform algebra $A$. 
Moreover, $\Delta_\prec+\Delta_\succ=\tdelta$ on this set of generators.
\item $\Delta_\prec$ and $\Delta_\succ$ are morphisms of $A$-dendriform modules.
\end{enumerate}
Then $A$ is a bidendriform bialgebra.
\end{prop}

{\bf Proof.}
Let us first show relation (\ref{E1'}). We put
$X=\{a\in A\:/\: \mbox{$a$ satisfies (\ref{E1'}))}\}$.
We have:
\begin{eqnarray*}
X&=& Ker\left((\Delta_\prec \otimes Id)\circ \Delta_\prec-
(Id \otimes \Delta_\prec+Id \otimes \Delta_\succ)\circ \Delta_\prec\right)\\
&=& Ker\left((\overline{\Delta}_\prec \overline{\otimes} Id)\circ \overline{\Delta}_\prec-
(Id \overline{\otimes} \overline{\Delta}_\prec+Id \overline{\otimes} \overline{\Delta}_\succ)
\circ \overline{\Delta}_\prec\right),
\end{eqnarray*}
which is the kernel of a certain morphism of dendriform modules from $A$
into $A\totimes A \totimes A$ by hypothesis 2. So $X$ is a dendriform submodule of $A$,
hence a dendriform subalgebra of $A$. As it contains a set of generators by hypothesis $1$, it is $A$.
We can prove (\ref{E2'}), (\ref{E3'}), and the fact that $\tdelta=\Delta_\prec+\Delta_\succ$
on the whole $A$ in the same way.
As the hypothesis 2 is a reformulation of axioms (\ref{E5})-(\ref{E8}),
$A$ is a bidendriform bialgebra. $\Box$

\subsection{Bidendriform structure on $\A^\D$}

Unfortunately, $(\A^\D,\prec,\succ, \Delta_\prec',\Delta_\succ')$
is not a bidendriform bialgebra. For example, for $a=b=\tdun{$d$}$:
$$
a\prec b=\tddeux{$d$}{$d$},\hspace{1cm}
\Delta'_\prec(a\prec b)=0,$$
whereas $a'b'_\prec\otimes a''\prec b''_\prec
+a'b \otimes a''+b'_\prec\otimes a\prec b''_\prec+b\otimes a
=\tdun{$d$} \otimes \tdun{$d$}$.
So (\ref{E8}) is not satisfied.

\begin{theo}
\label{theoADdend}
There is a unique structure of bidendriform bialgebra
(with the already known dendriform bialgebra structure) 
on $\A^\D$ such that for all $d\in \D$, 
$ \Delta_\prec(\tdun{$d$})= \Delta_\succ(\tdun{$d$})=0$. 
Hence, $(\A^\D,\prec,\succ,\Delta_\prec,\Delta_\succ)$ is a bidendriform bialgebra,
which induces the structure of Hopf algebra of $\h^\D$ already described.
\end{theo}

{\bf Proof.} We use the notations of proposition \ref{propo30}.

{\it Unicity.} As the $\tdun{$d$}$'s generate $\A^\D$,
there is at most one way to extend $\overline{\Delta}_\prec$
and $\overline{\Delta}_\succ$ to $A$ as morphisms of dendriform modules, with notations
of proposition \ref{propo30}.

{\it Existence.} By proposition \ref{propo29}, we can extend
$\overline{\Delta}_\prec$
and $\overline{\Delta}_\succ$ to $A$ as morphisms of dendriform modules.
So the second condition of proposition \ref{propo30} is satisfied.
The first one is trivially satisfied on the set of generators $\{\tdun{$d$}\:/\:d\in \D\}$. $\Box$\\

When $\D$ is a graded set, this structure obviously respects the gradation of $\A^\D$:

\begin{cor}
If $\D$ is a graded set, then $\A^\D$ is a graded bidendriform bialgebra.
\end{cor}

\begin{prop}
\label{propcoupe}
Let $F\in \F^D$, $F\neq 1$. We consider the set ${\cal A}dm_\prec(F)$ 
of admissible cuts of $F$ satisfying one of these two conditions:
\begin{enumerate}
\item $c$ cuts one of the edges which are on the path from the root
of the last tree of $F$ to the leave which is at most on the right of $F$.
\item $c$ cuts totally the last tree of $F$ if $F$ is not a single tree.
\end{enumerate}
Then:
\begin{eqnarray}
\label{eqdelta}
\Delta_\prec(F)&=&\sum_{c\in {\cal A}dm_\prec(F)}P^c(F)\otimes R^c(F).
\end{eqnarray}
\end{prop}

{\bf Prove.} We denote by $F'_\ll \otimes F''_\ll$ the second member of (\ref{eqdelta}).
Let us prove the result by induction on $n=weight(F)$. If $n=1$, then $F=\tdun{$d$}$ and the 
result is obvious. If $n\geq 2$, we have two possible cases:
\begin{enumerate}
\item $F=GH$, with $weight(G),weight(H)<n$. By a study of $\A dm_\prec(F)$, we easily have:
\begin{eqnarray*}
F_\ll'\otimes F_\ll'' &=&G'H'_\ll \otimes G''H''_\ll + G'H\otimes G''
+GH'_\ll \otimes H''_\ll+GH'_\ll \otimes H''_\ll+H\otimes G\\
 &=&G'H'_\prec \otimes G''H''_\prec + G'H\otimes G''
+GH'_\prec \otimes H''_\prec+GH'_\prec \otimes H''_\prec+H\otimes G,
\end{eqnarray*}
by the induction hypothesis on $H$. By (\ref{E78}), this is equal to $\Delta_\prec(F)$.
\item $F=B_d^+(G)$, $weight(G)=n-1$. By a study of $\A dm_\prec(G)$, we easily have:
\begin{eqnarray*}
F'_\ll\otimes F''_\ll &=&G'_\ll \otimes B_d^+(G''_\ll) +G\otimes \tdun{$d$}\\
&=&G'_\prec \otimes B_d^+(G''_\prec) +G\otimes \tdun{$d$}\\
&=&G'_\prec \otimes \tdun{$d$} \prec G''_\prec +G\otimes \tdun{$d$},
\end{eqnarray*}
by the induction hypothesis on $H$. As $F=\tdun{$d$}\prec G$, by (\ref{E8}) for $a=\tdun{$d$}$,
this is equal to $\Delta_\prec(F)$. $\Box$\\
\end{enumerate}

{\bf Example.}

\begin{eqnarray*}
\Delta_\prec \left( \tdun{$a$}\tdquatrequatre{$b$}{$c$}{$d$}{$e$}\right)&=&
 \tdquatrequatre{$b$}{$c$}{$d$}{$e$}\otimes \tdun{$a$}
+\tdun{$a$}\tdun{$e$}\tdun{$d$} \otimes \tddeux{$b$}{$c$}
+\tdun{$e$}\tdun{$d$} \otimes \tdun{$a$}\tddeux{$b$}{$c$}\\
&&+\tdun{$a$}\tdtroisun{$c$}{$d$}{$e$}\otimes \tdun{$b$}
+\tdtroisun{$c$}{$d$}{$e$}\otimes \tdun{$a$}\tdun{$b$}
+\tdun{$a$}\tdun{$d$} \otimes \tdtroisdeux{$b$}{$c$}{$e$}
+\tdun{$d$} \otimes \tdun{$a$}\tdtroisdeux{$b$}{$c$}{$e$},\\
\Delta_\succ \left( \tdun{$a$}\tdquatrequatre{$b$}{$c$}{$d$}{$e$}\right)&=&
\tdun{$a$}\otimes  \tdquatrequatre{$b$}{$c$}{$d$}{$e$}
+\tdun{$a$}\tdun{$e$} \otimes \tdtroisdeux{$b$}{$c$}{$d$}
+\tdun{$e$} \otimes\tdun{$a$} \tdtroisdeux{$b$}{$c$}{$d$}.
\end{eqnarray*}

\subsection{Totally primitive elements of $\A^\D$ and universal property}

We here prove that a connected bidendriform bialgebra is isomorphic to $\A^\D$
for a certain set $\D$.

\begin{prop}
\label{prop22}
 $Prim_{tot}(\A^\D)= Vect(\tdun{$d$}\:/\:d\in \D)$.
\end{prop}

{\bf Proof.} In a immediate way, $\tdun{$d$} \in Prim_{tot}(\A^\D)$.
Let us show the other inclusion.\\

{\it First case.}  Suppose that $\D$ is finite.
We graduate $\D$ by putting all its elements in degree $1$.
Then $\A^\D$ is graded, with $\A^\D_0=(0)$ and $\A^\D_n$ finite-dimensional for all $n$.
By corollary \ref{cor15}, $\A^\D=Prim_{tot}(\A^\D)\oplus(\A^\D)^{D2}$.
As $\A^\D$ is freely generated by the $\tdun{$d$}$'s, 
$Vect(\tdun{$d$}\:/\:d\in \D)\oplus (\A^\D)^{D2}=\A^\D$,
and this implies that $Prim_{tot}(\A^\D) \subseteq Vect(\tdun{$d$}\:/\:d\in \D)$. \\

{\it General case.} Let $p\in Prim_{tot}(\A^\D)$. There exists a finite subset $\D'$ of $\D$,
such that $p\in \A^{\D'}$. Then $p\in Vect(\tdun{$d$}\:/\:d\in \D')\subseteq 
Vect(\tdun{$d$}\:/\:d\in \D)$. $\Box$

\begin{theo}
\label{theo23}
Let $A$ be a bidendriform bialgebra. For all $d\in \D$,
let $p_d\in Prim_{tot}(A)$. There exists a unique morphism of bidendriform bialgebras:
$$ \Psi: \left\{\begin{array}{rcl}
\A^\D& \longrightarrow &A\\
\tdun{$d$}& \longrightarrow &p_d.
\end{array}\right. $$
Moreover, we have:
\begin{enumerate}
\item If the $p_d$'s are linearly independant, then $\Psi$ is monic.
\item if the dendriform coalgebra $A$ is connected and if the family $(p_d)_{d\in \D}$
linearly generates $Prim_{tot}(A)$, then $\Psi$ is epic. 
\item if the dendriform coalgebra $A$ is connected and if the family $(p_d)_{d\in \D}$
is a linear basis of $Prim_{tot}(A)$, then $\Psi$ is an isomorphism. 
\end{enumerate}
\end{theo}

{\bf Proof.} As $\A^\D$ is freely generated by the $\tdun{$d$}$'s, $\Psi$
defines a unique morphism of dendriform algebras.
As $\tdun{$d$}$ and $p_d$ are both totally primitive, $\Psi$ is a morphism of bidendriform
bialgebras. We now prove the two assertions.
\begin{enumerate}
\item We graduate $\D$ by putting all its elements of degree $1$. Suppose that $Ker(\Psi)\neq(0)$.
Let $x\in Ker(\Psi)$, non-zero, of minimal degree $n$. So $\Psi$ is monic on 
$\A^\D_{<n}=\A^\D_0\oplus\ldots\oplus \A^\D_{n-1}$. Moreover,
$\Delta_\prec(x)\in \A^\D_{<n}\otimes \A^\D_{<n}$ and
$(\Psi\otimes \Psi)\circ \Delta_\prec(x)=\Delta_\prec(\Psi(x))=0$.
By injectivity, $\Delta_\prec(x)=0$. In the same way, $\Delta_\succ(x)=0$.
So $x\in Prim_{tot}(\A^\D)=Vect(\tdun{$d$}\:/\:d\in \D)$ (proposition \ref{prop22}).
As the $p_d$'s are linearly independant, $Ker(\Psi)\cap Vect(\tdun{$d$}\:/\:d\in \D)=(0)$: 
contradiction. So $\Psi$ is monic.
\item Then $Im(\Psi)$ is a dendriform subalgebra of $A$ which contains $Prim_{tot}(A)$.
As $A$ is connected, $Prim_{tot}(A)$ generates $A$ (theorem \ref{theo17}), so $\Psi$ is epic.
\item Comes from 1 and 2. $\Box$
\end{enumerate} 

This corollary is immediate:
\begin{cor}
\label{cor24}
\begin{enumerate}
\item Let $A$ be a connected bidendriform bialgebra. Let $(p_d)_{d\in \D}$
be a basis of $Prim_{tot}(A)$. Then the morphism $\Psi:\A^\D \longrightarrow A$
described in theorem \ref{theo23} is an isomorphism of bidendriform bialgebras.
\item Let $A$ be a $\mathbb{N}$-graded bidendriform bialgebra which is connected
as a dendriform coalgebra. Let $(p_d)_{d\in \D}$
be a basis of $Prim_{tot}(A)$ made of homogeneous elements. 
Then $\D$ is given a gradation by putting $|d|=|p_d|$.
Then the morphism $\Psi:\A^\D \longrightarrow A$
described in theorem \ref{theo23} is an isomorphism of graded bidendriform bialgebras.
\end{enumerate}
\end{cor}

\begin{cor}
\label{cor25}
Let $A$ be a $\mathbb{N}$-graded bidendriform bialgebra such that $A_0=(0)$
and $A_n$ is finite-dimensional for all $n\in \mathbb{N}$.
We consider the following formal series:
$$P(X)=\sum_{n=1}^{+\infty}dim(Prim_{tot}(A)_n)X^n,\hspace{1cm}
R(X)=\sum_{n=1}^{+\infty}dim(A_n)X^n.$$
Then $\displaystyle P(X)=\frac{R(X)}{(1+R(X))^2}$.
\end{cor}

{\bf Proof.} 
Immediate if $A=(0)$. Suppose that $A\neq(0)$ (so $P(X)\neq 0$).
By corollary \ref{cor24}, we can suppose $A=\A^\D$. 
By proposition \ref{prop22},
$P(X)=D(X)$ of proposition \ref{prop18}. Hence: 
\begin{eqnarray*}
&& R(X)+1=\frac{1-\sqrt{1-4P(X)}}{2P(X)}\\
&\Rightarrow&  2P(X)(R(X)+1)-1=-\sqrt{1-4P(X)}\\
&\Rightarrow&  4P(X)^2(R(X)+1)^2+1-4P(X)(R(X)+1)=1-4P(X)\\
&\Rightarrow&  P(X)(R(X)+1)^2-R(X)-1=-1\\
&\Rightarrow&  P(X)=\frac{R(X)}{(R(X)+1)^2}.\:\Box
\end{eqnarray*}

\section{Application to the Hopf algebra $\malv$}

\label{sect4} 

\subsection{Recalls}

(See \cite{Aguiar2,Duchamp,Malvenuto}).
The algebra $\malv$ is the vector space generated by the 
elements $(\FF_u)_{u\in \S}$, where $\S$
is the disjoint union of the symmetric groups $S_n$ $(n\in \mathbb{N})$.
Its product and its coproduct are given in the following way:
for all $u\in S_n$, $v\in S_m$, by putting $u=(u_1\ldots u_n)$,
\begin{eqnarray*}
\Delta(\FF_u)&=&\sum_{i=0}^n \FF_{st(u_1\ldots u_i)}\otimes \FF_{st(u_{i+1}\ldots u_n)},\\
\FF_u.\FF_v&=&\sum_{\zeta\in sh(n,m)}\FF_{(u\times v).\zeta^{-1}},
\end{eqnarray*}
where $sh(n,m)$ is the set of $(n,m)$-shuffles, and $st$ is the standardisation.
Its unit is $1=\FF_\emptyset$, where $\emptyset$ is the unique element of $S_0$.
Moreover, $\malv$ is a $\mathbb{N}$-graded Hopf algebra, by putting $|\FF_u|=n$
if $u\in S_n$.\\

{\bf Examples.}
\begin{eqnarray*}
\FF_{(1\:2)}\FF_{(1\:2\:3)}&=& 
\FF_{(1\:2\:3\:4\:5)}+ \FF_{(1\:3\:2\:4\:5)}+ \FF_{(1\:3\:4\:2\:5)}+ \FF_{(1\:3\:4\:5\:2)}
+\FF_{(3\:1\:2\:4\:5)}\\
&&+ \FF_{(3\:1\:4\:2\:5)}+ \FF_{(3\:1\:4\:5\:2)}+\FF_{(3\:4\:1\:2\:5)}
+\FF_{(3\:4\:1\:5\:2)} +\FF_{(3\:4\:5\:1\:2)},\\
\\
\Delta\left(\FF_{(1\:2\:5\:4\:3)}\right)&=& 
1\otimes \FF_{(1\:2\:5\:4\:3)}+\FF_{(1)} \otimes \FF_{(1\:4\:3\:2)}
+\FF_{(1\:2)} \otimes \FF_{(3\:2\:1)}\\
&&+\FF_{(1\:2\: 3)} \otimes \FF_{(2\:1)}
+\FF_{(1\:2\:4\:3)} \otimes \FF_{(1)}+ \FF_{(1\:2\:5\:4\:3)}\otimes 1.
\end{eqnarray*}

\subsection{Bidendriform structure on $\malv$}

Let $(\malv)_+=Vect(\FF_u\:/\:u\in S_n,\: n\geq 1)$ be the augmentation ideal of $\malv$. 
We define $\prec,\succ,\Delta_\prec$ and $\Delta_\succ$
on $(\malv)_+$ in the following way: for all $u\in S_n$, $v\in S_m$,
by putting $u=(u_1\ldots u_n)$,
\begin{eqnarray*}
\FF_u\prec\FF_v&=&\sum_{\stackrel{\zeta\in sh(n,m)}{\zeta^{-1}(n+m)=n}}\FF_{(u\times v).\zeta^{-1}},\\
\FF_u\succ\FF_v&=&\sum_{\stackrel{\zeta\in sh(n,m)}{\zeta^{-1}(n+m)=n+m}}\FF_{(u\times v).\zeta^{-1}},\\
\Delta_\prec(\FF_u)&=&\sum_{i=u^{-1}(n)}^{n-1} \FF_{st(u_1\ldots u_i)}\otimes \FF_{st(u_{i+1}\ldots u_n)},\\
\Delta_\succ(\FF_u)&=&\sum_{i=1}^{u^{-1}(n)-1} \FF_{st(u_1\ldots u_i)}\otimes \FF_{st(u_{i+1}\ldots u_n)}.
\end{eqnarray*}

{\bf Examples.}
\begin{eqnarray*}
\FF_{(1\:2)}\prec \FF_{(1\:2\:3)}&=& \FF_{(1\:3\:4\:5\:2)} +\FF_{(3\:1\:4\:5\:2)} 
+\FF_{(3\:4\:1\:5\:2)} +\FF_{(3\:4\:5\:1\:2)},\\
\FF_{(1\:2)}\succ \FF_{(1\:2\:3)}&=& \FF_{(1\:2\:3\:4\:5)} +\FF_{(1\:3\:2\:4\:5)} 
+\FF_{(1\:3\:4\:2\:5)}+\FF_{(3\:1\:2\:4\:5)} +\FF_{(3\:1\:4\:2\:5)} +\FF_{(3\:4\:1\:2\:5)}, \\
\\
\Delta_\prec\left(\FF_{(1\:2\:5\:4\:3)}\right)&=& \FF_{(1\:2\: 3)} \otimes \FF_{(2\:1)}
+\FF_{(1\:2\:4\:3)} \otimes \FF_{(1)},\\
\Delta_\succ\left(\FF_{(1\:2\:5\:4\:3)}\right)&=& \FF_{(1)} \otimes \FF_{(1\:4\:3\:2)}
+\FF_{(1\:2)} \otimes \FF_{(3\:2\:1)}.
\end{eqnarray*}

\begin{theo}
($(\malv)_+,\prec,\succ,\Delta_\prec,\Delta_\succ)$ is a connected bidendriform bialgebra.
\end{theo}

{\bf Proof.} The structure of dendriform bialgebra is already introduced
in \cite{Ronco2}, so we already have (\ref{E1})-(\ref{E3}), and (\ref{E5})+(\ref{E7}),
(\ref{E6})+(\ref{E8}).
We consider the symmetric non-degenerate pairing on $(\malv)_+$ defined in \cite{Duchamp,Malvenuto} by
$<\FF_\sigma,\FF_\tau>=\delta_{\sigma,\tau^{-1}}$.
Note that $\Delta_\prec,\Delta_\succ$ are the transposes of $\prec,\succ$ for this pairing.
So $((\malv)_+,\Delta_\prec,\Delta_\succ)$ is a codendriform coalgebra.
So we already have (\ref{E1'})-(\ref{E3'}), and (\ref{E5})+(\ref{E6}),
(\ref{E7})+(\ref{E8}). It is then enough to prove (\ref{E6}).
We have:
\begin{eqnarray*}
\Delta_\succ\left(\FF_{(u_1\ldots u_n)}\prec \FF_{(v_1\ldots v_m)}\right)
&=&\Delta_\succ\left(\sum_{1\leq j_1< \ldots <j_m<n+m}
\FF_{\underbrace{(u_1\ldots v_1+n\ldots  v_m+n \ldots u_n)}_{\mbox{$v_i+n$ in position $j_i$}}}\right)\\
&=&\sum \sum_{1\leq j_1< \ldots <j_m<n+m}
\FF_{st(u_1\ldots)}\otimes \FF_{st(\ldots m+n\ldots u_n)}\\
&=&\underbrace{\left(\FF_{(u_1\ldots u_n)}\right)' \left(\FF_{(v_1\ldots v_n)}\right)'_\succ
\otimes \left(\FF_{(u_1\ldots u_n)}\right)''\prec \left(\FF_{(v_1\ldots v_n)}\right)''_\succ}_
{\mbox{terms with $u_i$'s and $v_j$'s on both sides on the $\otimes$}}\\
&&+\underbrace{\left(\FF_{(u_1\ldots u_n)}\right)' 
\otimes \left(\FF_{(u_1\ldots u_n)}\right)''\prec \FF_{(v_1\ldots v_n)}}_
{\mbox{terms with all the $v_i$'s on the right on the $\otimes$}}\\
&+&\underbrace{ \left(\FF_{(v_1\ldots v_n)}\right)'_\succ
\otimes \FF_{(u_1\ldots u_n)}\prec \left(\FF_{(v_1\ldots v_n)}\right)''_\succ}_
{\mbox{terms with all the $u_i$'s the right of the $\otimes$}}.
\end{eqnarray*}
So (\ref{E6}) is satisfied, and $(\malv)_+$ is a bidendriform bialgebra. Moreover,
it is $\mathbb{N}$-graded, with $((\malv)_+)_0=(0)$, so it is connected. $\Box$\\

{\bf Remark.} It is of course possible to prove (\ref{E1'})-(\ref{E3'}) directly. 
For example, the two members of (\ref{E1'}) for $a=\FF_u$, with $u=(u_1\ldots u_n)$,
are equal to:
$$\sum_{\sigma^{-1}(n)\leq i<j<n}
\FF_{st(u_1\ldots u_i)}\otimes \FF_{st(u_{i+1}\ldots u_j)} \otimes \FF_{st(u_{j+1}\ldots u_n)}.$$
For (\ref{E2'}), we obtain:
$$\sum_{1\leq i<\sigma^{-1}(n)\leq j<n}
\FF_{st(u_1\ldots u_i)}\otimes \FF_{st(u_{i+1}\ldots u_j)} \otimes \FF_{st(u_{j+1}\ldots u_n)}.$$
And for (\ref{E3'}):
$$\sum_{1\leq i<j<\sigma^{-1}(n)}
\FF_{st(u_1\ldots u_i)}\otimes \FF_{st(u_{i+1}\ldots u_j)} \otimes \FF_{st(u_{j+1}\ldots u_n)}.$$
\\

Moreover, $(\malv)_+$ is graded and connected, with:
$$R(X)=\sum_{n=1}^{+\infty} dim(\malv_n)X^n=\sum_{n=1}^{+\infty}n!X^n.$$
By corollary \ref{cor25}:
$$P(X)=\sum_{n=1}^{+\infty}dim(Prim_{tot}((\malv)_+)_n)X^n=\sum_{n=1}^{+\infty}p_nX^n
=\frac{R(X)}{(R(X)+1)^2}.$$
We obtain then:

$$\begin{array}{|c|c|c|c|c|c|c|c|c|c|c|c|c|}
\hline
n&1&2&3&4&5&6&7&8&9&10&11&12\\
\hline
p_n&1&0&1&6&39&284&2\:305&20\:682&203\:651&2\:186\:744&25\: 463 \:925
&319\: 989\: 030\\
\hline
\end{array}$$

{\bf Remark.} Let $\Psi:\A\longrightarrow (\malv)_+$ be the unique morphism of bidendriform bialgebras
which send $\tun$ to $\FF_{(1)}$. It coincides with the morphism $\overline{\Phi^*}$
of \cite{Loday3}. We obtain by theorem \ref{theo23} that it is monic.\\

By corollary \ref{cor24}:
\begin{theo}
\label{theo39}
Let $\D$ be a graded set such that $\D_0=\emptyset$ and for all $n\geq 1$,
$card(\D_n)=p_n$. Then $(\malv)_+$ and $\A^\D$ are isomorphic as  graded bidendriform bialgebras.
Hence, $\malv$ and $\h^\D$ are isomorphic as graded Hopf algebras.
\end{theo}

We can now prove the conjecture 3.8 of \cite{Duchamp}:
\begin{cor}
\label{cor40}
If the field $K$ is of characteristic $0$, the Lie algebra $Prim(\malv)$ is free.
\end{cor}

{\bf Proof.} By proposition 141 of \cite{Foissy4},
the Lie algebra $Prim(\h^\D)$ is free. $\Box$\\

Finally, we give a basis of totally primitive elements of $(\malv)_+$ of degree smaller than $4$.
With notations of \cite{Duchamp}:
\begin{eqnarray*}
\V_1&=&\FF_1,\\
\\
\V_{231}&=&\FF_{231}-\FF_{132},\\
\\
\V_{3142}&=&\FF_{3142}-\FF_{2143},\\
\V_{2431}&=&\FF_{2431}-\FF_{1432},\\
\V_{2341}&=&\FF_{2341}-\FF_{1342},\\
\V_{3241}&=& \FF_{1243}-\FF_{1342}-\FF_{2143}+\FF_{3241},\\
\V_{3412}-\V_{2413}&=&\FF_{3412}-\FF_{2413},\\
\V_{3421}-\V_{2413}&=&\FF_{1423}-\FF_{2413}-\FF_{1432}+\FF_{3421}.
\end{eqnarray*}

\section{Bidendriform pairings}

\subsection{Definition}

\begin{defi}
\textnormal{Let $A$ and $B$ be two bidendriform  bialgebras
and let $<,>:A\times B \longrightarrow K$ be a bilinear form.
We will say that $<,>$ is a bidendriform pairing if, for all $a,a_1,a_2\in A$, $b,b_1,b_2\in B$,
the following conditions are satisfied:
$$\begin{array}{rclc|crcl}
<a_1\prec a_2,b>&=&<a_1\otimes a_2,\Delta_\prec(b)>,&&&
<a,b_1\prec b_2>&=&<\Delta_\prec(a),b_1\otimes b_2>,\\
<a_1\succ a_2,b>&=&<a_1\otimes a_2,\Delta_\succ(b)>,&&&
<a,b_1\succ b_2>&=&<\Delta_\succ(a),b_1\otimes b_2>.
\end{array}$$
Suppose that $A$ and $B$ are graded. Then $<,>$ will be said homogeneous if for all
homogeneous elements  $a\in A$, $b\in B$, $<a,b>=0$ if $a$ and $b$ have different degrees.}
\end{defi}

{\bf Remark.} If $<,>:A\times B \longrightarrow K$ is a bidendriform pairing,
then it induces a Hopf pairing between $\overline{A}=A\oplus K$ and $\overline{B}=B\oplus K$ 
with, for all $a\in A$, $b\in B$,
$<1,b>=0$, $<a,1>=0$, and $<1,1>=1$.\\

{\bf Example.} The pairing on $(\malv)_+$ defined by $<\FF_u,\FF_v>=\delta_{u,v^{-1}}$
is a bidendriform pairing.

\subsection{Hopf dual of a bidendriform bialgebra}

We here extend the notion of Hopf dual of a Hopf algebra (see \cite{Abe,Joseph}) to bidendriform bialgebras.
Let $A$ be a dendriform algebra. 
We denote by $A^\star$ the following subspace of the linear dual $A^*$ of $A$:
$$A^\star=\{f\in A^*\:/\: \mbox{$f$ vanishes on a dendriform ideal of finite codimension of $A$}\}.$$

For all $a\in A$, we consider the following applications:

$$ L_a^\prec: \left\{\begin{array}{rcl}
A^* & \longrightarrow &A^*\\
f& \longrightarrow & 
\left\{\begin{array}{rcl}
A& \longrightarrow &K\\
b &\longrightarrow &L_a^\prec(f)(b)=f(a\prec b),
\end{array}\right. 
\end{array}\right. $$
$$ L_a^\succ: \left\{\begin{array}{rcl}
A^* & \longrightarrow &A^*\\
f& \longrightarrow & 
\left\{\begin{array}{rcl}
A& \longrightarrow &K\\
b &\longrightarrow &L_a^\succ(f)(b)=f(a\succ b),
\end{array}\right. 
\end{array}\right. $$
$$ R_a^\prec: \left\{\begin{array}{rcl}
A^* & \longrightarrow &A^*\\
f& \longrightarrow & 
\left\{\begin{array}{rcl}
A& \longrightarrow &K\\
b &\longrightarrow &R_a^\prec(f)(b)=f(b\prec a),
\end{array}\right. 
\end{array}\right. $$
$$ R_a^\succ: \left\{\begin{array}{rcl}
A^* & \longrightarrow &A^*\\
f& \longrightarrow & 
\left\{\begin{array}{rcl}
A& \longrightarrow &K\\
b &\longrightarrow &R_a^\succ(f)(b)=f(b\succ a).
\end{array}\right. 
\end{array}\right. $$

\begin{lemme}
\label{lemmedual}
Let $f\in A^*$. The following assertions are equivalent:
\begin{enumerate}
\item $f \in A^\star$.
\item There exists a finite dimensional subspace $U$ of $A^*$, stable for all $a\in A$
by $L_a^\prec$, $L_a^\succ$, $R_a^\prec$ and $R_a^\succ$, which contains $f$.
\end{enumerate}
\end{lemme}

{\bf Proof.} 

$1\Longrightarrow 2$. Let $I$ be a dendriform ideal of $A$, such that $f(I)=(0)$ and $f$ has
a finite codimension. Then $I^\perp$ is finite-dimensional and contains $f$. 
Let $g\in I^\perp$ and let $a\in A$. For all $b\in I$:
$$L_a^\prec(g)(b)=g(\underbrace{a\prec b}_{\in I})=0,$$
as $I$ is a dendriform ideal of $A$. So $L_a^\prec(g)\in I^\perp$: $I^\perp$ is stable under $L_a^\prec$.
In the same way, $I^\perp$ is stable under $L_a^\succ$, $R_a^\prec$ and $R_a^\succ$.\\

$2 \Longrightarrow 1$. Let $I=U^\perp$. Then $I$ has a finite codimension. 
Moreover, as $f \in U$, $f(I)=(0)$. Let $a\in A$, $b\in I$. For all $g\in U$:
$$g(a\prec b)=\underbrace{L_a^\prec(g)}_{\in U}(b)=0,$$
so $a\prec b \in I$. In the same way, $a\succ b$, $b\prec a$ and $b\succ a \in I$, 
so $I$ is a dendriform ideal. $\Box$\\

\begin{lemme}
$\prec^*$, $\succ^*$ send $A^\star$ into $A^\star \otimes A^\star$.
\end{lemme}

{\bf Proof.} Let $f\in A^\star$ and $U$ as in lemma \ref{lemmedual}.
Let $(u_1,\ldots,u_n)$ be a basis of $U$. By lemma \ref{lemmedual}, observe that $U \subseteq A^\star$.
For all $a\in A$, we put:
$$L_a^\prec(f)=\sum_{i=1}^n v_i(a)u_i.$$
Then the $v_i$'s are elements of $A^*$.
Let us fix $j\in \{1,\ldots,n\}$. There exists $b_j \in A$, such that $v_i(b_j)=\delta_{i,j}$. Then:
$$L_a^\prec(f)(b_j)=\sum_{i=1}^n v_i(a)u_i(b_j)=v_j(a)=f(a\prec b_j)=R_{b_j}^\prec(f)(a),$$
so $v_j=R_{b_j}^\prec(f) \in U$. Moreover, for all $a,b\in A$:
$$\prec^*(f)(a\otimes b)=f(a\prec b)=L_a^\prec(f)(b)=\sum_{i=1}^n v_i(a)\otimes u_i(b)
=\sum_{i=1}^n v_i \otimes u_i (a\otimes b),$$
so $\prec^*(f) \in U\otimes U\subseteq A^\star \otimes A^\star$. In the same way,
$\succ^*(f) \in A^\star \otimes A^\star$. $\Box$\\

{\bf Remark.} So $(A^\star,\prec^*,\succ^*)$ is a dendriform coalgebra.\\

Suppose now that $A$ is a bidendriform bialgebra.
By duality, $(A^*,\Delta_\prec^*,\Delta_\succ^*)$ is given a dendriform algebra structure.
More precisely, for all $f,g\in A^*$, we have, for all $a\in A$:
$$(f\prec g)(a)=(f\otimes g)\circ\Delta_\prec(a),\hspace{1cm}
(f\succ g)(a)=(f\otimes g)\circ\Delta_\succ(a).$$

\begin{lemme}
$A^\star$ is a dendriform subalgebra of $A^*$. 
\end{lemme}

{\bf Proof.} Let $f,g\in A^\star$.
By lemma \ref{lemmedual}, there exists a finite-dimensional subspace $U$ of $A^*$,
containing $f$ and $g$, stable for all $a\in A$
by $L_a^\prec$, $L_a^\succ$, $R_a^\prec$ and $R_a^\succ$. 
Let $F,G \in U$. Let $a\in A$. For all $b\in A$, by (\ref{E8}):
\begin{eqnarray*}
L_a^\prec(F\prec G)(b)&=&(F\prec G)(a\prec b)\\
&=&(F\otimes G)\circ \Delta_\prec(a\prec b)\\
&=&(F\otimes G)(a'b'_\prec \otimes a''\prec b''_\prec
+a'b\otimes a''+b'_\prec \otimes a\prec b''_\prec+b\otimes a)\\
&=&\left((L_{a'}^\prec+L_{a'}^\succ)(F)\prec L_{a''}^\prec(G)\right)(b)\\
&&+G(a'')\left((L_{a'}^\prec+L_{a'}^\succ)(F)\right)(b)+ \left(F\prec L_a^\prec(G)\right)(b)+F(a)G(b).
\end{eqnarray*}
Hence:
\begin{eqnarray*}
L_a^\prec(F\prec G)&=&(L_{a'}^\prec+L_{a'}^\succ)(F)\prec L_{a''}^\prec(G)
+G(a'')(L_{a'}^\prec+L_{a'}^\succ)(F)+ F\prec L_a^\prec(G)+F(a)G\\
&\in&U+U\prec U+U\succ U.
\end{eqnarray*}
In the same way, we can prove that $V=U+U\prec U+U\succ U$ is stable under $L_a^\succ$,
$R_a^\prec$ and $R_a^\succ$. As $U$ is finite-dimensional, $V$ is finite-dimensional,
and contains $f\prec g$ and $f\succ g$, as $f,g\in U$. So, by lemma \ref{lemmedual},
$f\prec g$ and $f\succ g \in A^\star$. $\Box$\\

By dualising the axioms of  bidendriform bialgebras, we obtain:

\begin{prop}
For any bidendriform bialgebra $A$, $(A^\star, \Delta_\prec^*,\Delta_\succ^*,
\prec^*,\succ^*)$ is a bidendriform bialgebra. More precisely, for all $f,g\in A^\star$, $a,b\in A$:
$$(f\prec g)(a)=(f\otimes g)\circ\Delta_\prec(a),\hspace{1cm}
(f\succ g)(a)=(f\otimes g)\circ\Delta_\succ(a),$$
$$\Delta_\prec(f)(a\otimes b)=f(a\prec b),\hspace{1cm}
\Delta_\succ(f)(a\otimes b)=f(a\succ b).$$
\end{prop}

\begin{prop}
\label{propo46}
Let $A,B$ be bidendriform bialgebras.
\begin{enumerate}
\item Let $\Phi:A\longrightarrow B^\star$ be a morphism of bidendriform bialgebras. Then the following pairing
is a bidendriform pairing:
$$<,>: \left\{\begin{array}{rcl}
A\times B& \longrightarrow &K\\
(a,b)& \longrightarrow &<a,b>=\Phi(a)(b)
\end{array}\right. $$
\item Suppose that $A$ is connected and that $<,>:A\times B \longrightarrow K$ is a bidendriform pairing.
Then the following application is a morphism of bidendriform bialgebras:
$$\Phi: \left\{\begin{array}{rcl}
A& \longrightarrow &B^\star\\
a& \longrightarrow &\Phi(a):\left\{\begin{array}{rcl}
B& \longrightarrow &K\\
b& \longrightarrow &<a,b>.
\end{array}\right.
\end{array}\right. $$
\end{enumerate}
\end{prop}

{\bf Proof.}

1. Let $a_1,a_2\in A$, $b\in B$.
\begin{eqnarray*}
<a_1\prec a_2,b>&=&\Phi(a_1\prec a_2)(b)\\
&=&(\Phi(a_1)\prec \Phi(a_2))(b)\\
&=&(\Phi(a_1)\otimes \Phi(a_2))\circ \Delta_\prec(b)\\
&=&<a_1\otimes a_2,\Delta_\prec(b)>.
\end{eqnarray*}

Let $a\in A$, $b_1,b_2 \in B$.
\begin{eqnarray*}
<a,b_1\prec b_2>&=&\Phi(a)(b_1\prec b_2)\\
&=&\Delta_\prec(\Phi(a))(b_1 \otimes b_2)\\
&=&((\Phi \otimes \Phi)\circ \Delta_\prec(a))(b_1\otimes b_2)\\
&=&<\Delta_\prec(a),b_1\otimes b_2>.
\end{eqnarray*}
The two other equalities are proved in the same way.\\

2. We first prove the following lemma:

\begin{lemme}
\label{lem47}
Let $C$ be a connected dendriform coalgebra. For all $a\in C$, 
there exists a finite-dimensional dendriform subcoalgebra of $C$ which contains $a$.
\end{lemme}

{\bf Proof.} We use the notations of definition \ref{deficonnexe}.
Let $a\in C$. We put:
$$C_a=Vect\left\{\begin{array}{c}
(f_1\otimes \ldots \otimes f_{i-1}\otimes Id \otimes f_{i+1} \otimes \ldots \otimes f_n)
\circ P(a)\:/\\
 n\in \mathbb{N}^*, \:i\in \{1,\ldots,n\}, \:f_1,\ldots,f_{i-1},f_{i+1},\ldots f_n \in C^*, \:
P\in \P_C(n)\end{array}\right\}.$$
Then $C_a$ is the smallest dendriform subcoalgebra of $C$ which contains $a$.
As $C$ is connected, there are only a finite number of $P$ such that $P(a)\neq 0$,
so $C_a$ is finite-dimensional. $\Box$\\

{\bf End of the proof of proposition \ref{propo46}.}
Let us first prove that $\Phi(a)\in B^\star$ for all $a\in A$.
As $A$ is connected, there exists a finite-dimensional subcoalgebra $C$ of $A$ which contains $a$
(lemma \ref{lem47}). Let $I$ be the orthogonal of $C$ for the pairing $<,>$. As $C$ is finite-dimensional,
$I$ has a finite codimension. Let $b_1\in B$, $b_2 \in I$. For all $c\in C$, we have:
$$<c,b_1\prec b_2>=<\underbrace{\Delta_{\prec}(c)}_{\in C\otimes C},b_1\otimes b_2>=0,$$
so $b_1\prec b_2 \in I$. In the same way, we prove that $I$ is a dendriform ideal of $B$.
As $a\in C$, $\Phi(a)(I)=<a,I>=(0)$, so $\Phi(a)\in B^\star$.

Let $a_1,a_2 \in A$. For all $b\in B$:
\begin{eqnarray*}
\Phi(a_1\prec a_2)(b)&=&<a_1\prec a_2,b>\\
&=&<a_1\otimes a_2,\Delta_\prec(b)>\\
&=&(\Phi(a_1)\otimes \Phi(a_2))\circ\Delta_\prec(b)\\
&=&(\Phi(a_1)\prec \Phi(a_2))(b),
\end{eqnarray*}
so $\Phi(a_1\prec a_2)=\Phi(a_1)\prec \Phi(a_2)$. Let $a\in A$. For all $b_1,b_2\in B$:
\begin{eqnarray*}
(\Delta_\prec(\Phi)(a))(b_1 \otimes b_2)
&=&\Phi(a)(b_1 \prec b_2)\\
&=&<a,b_1\prec b_2>\\
&=&<\Delta_\prec(a),b_1 \otimes b_2>\\
&=&((\Phi \otimes \Phi)\circ\Delta_\prec(a))(b_1\otimes b_2),
\end{eqnarray*}
so $\Delta_\prec(\Phi)(a)=(\Phi \otimes \Phi)\circ\Delta_\prec(a)$.
The other equalities are proved in the same way. $\Box$

\subsection{Bidendriform pairing on a connected graded bidendriform bialgebra}

\begin{theo}
\label{theo30}
Let $A,B$ be connected bidendriform bialgebras.
Let $<,>_{|Prim_{tot}}:Prim_{tot}(A)\times Prim_{tot}(B)\longrightarrow K$. then
$<,>_{|Prim_{tot}}$ can be uniquely extended into a bidendriform pairing between $A$ and $B$.
Moreover:\begin{enumerate}
\item $<,>$ is non-degenerate $\Longleftrightarrow$ $<,>_{\mid Prim_{tot}}$ is non-degenerate.
\item Suppose that $A=B$. Then $<,>$ is symmetric $\Longleftrightarrow$ $<,>_{\mid Prim_{tot}}$ is symmetric.
\item Suppose that $A$ and $B$ are $\mathbb{N}$-graded
(we do not suppose that the $A_n$'s or the $B_n$'s are finite-dimensional). Then
$<,>$ is homogeneous $\Longleftrightarrow$ $<,>_{\mid Prim_{tot}}$ is homogeneous.
\end{enumerate}
\end{theo}

{\bf Proof.} 

{\it Existence.} 
By corollary \ref{cor15}, $B=Prim_{tot}(B)\oplus B^{D2}$.
We can then define, for all $p\in Prim_{tot}(A)$:
$$ f_p: \left\{\begin{array}{rcl}
B& \longrightarrow &K\\
q\in Prim_{tot}(B)& \longrightarrow &<p,q>_{\mid Prim_{tot}},\\
b\in B^{D2}&\longrightarrow 0.
\end{array}\right. $$
Then $Ker(f_p)$ contains $B^{D2}$, so is a dendriform ideal of codimension at most $1$.
Hence, $f_p \in B^\star$.
Moreover, for all $a_1,a_2 \in A$,
$\Delta_\prec(f_p)(a_1\otimes a_2)=f_p(a_2\prec a_2)=0$,
so $\Delta_\prec(f_p)=0$. In the same way, $\Delta_\succ(f_p)=0$ and then $f_p \in Prim_{tot}(B^\star)$.

As $A$ is isomorphic to a bidendriform algebra $\A^D$ (corollary \ref{cor24}),
by theorem \ref{theo23}, there exists a unique morphism
of bidendriform bialgebras:
$$ \Phi: \left\{\begin{array}{rcl}
A &\longrightarrow &B^\star\\
p\in Prim_{tot}(A)& \longrightarrow &f_p.
\end{array}\right. $$
We then put $<a,b>=\Phi(a)(b)$ for all $a\in A$, $b\in B$.
This pairing extends $<,>_{\mid Prim_{tot}}$: if $p \in Prim_{tot}(A)$, $q\in Prim_{tot}(B)$,
$<p,q>=\Phi(p)(q)=f_p(q)=<p,q>_{\mid Prim_{tot}}$.
Moreover, as $\Phi$ is a morphism of  bidendriform bialgebras,
$<,>$ is a bidendriform pairing (proposition \ref{propo46}).\\

{\it Unicity.} Let $<,>'$ a second bidendriform pairing which extends $<,>_{|Prim_{tot}}$.
Let us show that it is equal to $<,>$ already defined.
Consider:
$$A'=\{a\in A\:/\:\forall b\in B,\: <a,b>=<a,b>'\}.$$
Let $a_1,a_2 \in A'$. Then, forall $b\in B$,
$$<a_1\prec a_2,b>=<a_1\otimes a_2,\Delta_\prec(b)>=<a_1\otimes a_2,\Delta_\prec(b)>'=<a_1\prec a_2,b>'.$$
So $a_1\prec a_2\in A'$. In the same way, $a_1\succ a_2 \in A'$,
so $A'$ is a dendriform subalgebra.
As $A$ is generated by $Prim_{tot}(A)$, it is enough to prove that
$p\in A'$ for all $p \in Prim_{tot}(A)$. Let $b\in B$ and  let us show that $<p,b>=<p,b>'$.
It is true by hypothesis if $b\in Prim_{tot}(A)$.
If it is not the case, we can suppose $b=b_1\prec b_2$ or $b=b_1\succ b_2$.
Hence:
$$
\begin{array}{rcccccccl}
<p,b>'&=&<p,b_1\prec b_2>'&=&<\Delta_\prec(p),b_1\otimes b_2>'&=&0&=&<p,b>,\\
\mbox{or }<p,b>'&=&<p,b_1\succ b_2>'&=&<\Delta_\succ(p),b_1\otimes b_2>'&=&0&=&<p,b>.\\
\end{array}$$

We now prove the three equivalences.\\

{\it Non-degeneracy, $\Longrightarrow$.}
Suppose $<,>$ non-degenerate and let $a \in Prim_{tot}(A)\cap Prim_{tot}(B)^\perp$.
Then, as $Prim_{tot}(A)\subseteq (B^{D2})^\perp$,
we obtain $a\in (B^{D2}+Prim_{tot}(B))^\perp=B^\perp=(0)$.
In the same way, $Prim_{tot}(B) \cap Prim_{tot}(A)^\perp=(0)$,
So $<,>_{|Prim_{tot}}$ is  non-degenerate.

{\it Non-degeneracy, $\Longleftarrow$.} Then $\Phi_{\mid Prim_{tot}(A)}$ 
is injective, so, by theorem \ref{theo23}, $\Phi$ is injective: $B^\perp=(0)$.
By permuting the role of $A$ and $B$, which does not change the pairing $<,>$ by unicity,
we obtain $A^\perp=(0)$. So $<,>$ is  non-degenerate.\\

{\it Symmmetry, $\Longrightarrow$.} Obvious.

{\it Symmetry, $\Longleftarrow$.} We consider
$A''=\{a\in A\:/\:\forall b\in A,\: <a,b>=<b,a>\}$.
Let $a_1,a_2 \in A''$. Then, for all $b\in A$: 
$$<a_1\prec a_2,b>=<a_1\otimes a_2,\Delta_{\prec}(b)>
=<\Delta_{\prec}(b), a_1\otimes a_2>=<b,a_1\prec a_2>.$$
So $a_1\prec a_2 \in A''$.
In the same way,  $a_1\succ a_2 \in A''$, so $A''$ is a dendriform subalgebra.
It is then enough to show that $Prim_{tot}(A)\subset A''$.
Let $a\in Prim_{tot}(A)$ and $b\in A$. 
We can limit ourselves to the three following cases:
\begin{enumerate}
\item $b\in Prim_{tot}(A)$. Then $<a,b>=<b,a>$ by hypothesis on $<,>_{|Prim_{tot}}$.
\item $b=b_1\prec b_2$. Then $<b,a>=<b_1\otimes b_2,\Delta_{\prec}(a)>=0=<a,b>$.
\item $b=b_1\succ b_2$. Then $<b,a>=<b_1\otimes b_2,\Delta_{\succ}(a)>=0=<a,b>$.
\end{enumerate}
So $A''=A$ and $<,>$ is symmetric.\\

{\it Homogeneity, $\Longrightarrow$.} Obvious.

{\it Homogeneity, $\Longleftarrow$.} Then  $\Phi$ takes its values in the graded dual
of $A$ and is homogeneous of degree $0$, so $<,>$ is homogeneous. $\Box$

\subsection{Application to $\h^\D$}

We define a pairing on $Prim_{tot}(\A^\D)$ in the following way:
for all $d,d'\in \D$,
$$<\tdun{$d$},\tdun{$d'$}>_{Prim_{tot}(\A^\D)}=\delta_{d,d'}.$$
This pairing is symmetric, homogeneous, and non-degenerate.
By theorem \ref{theo30}, it can be extended on a unique bidendriform pairing $<,>$ on $\A^\D$.
This pairing is also symmetric, homogeneous, and non-degenerate.\\

{\bf Remark.} This pairing is not the pairing $(,)$ on $\h^\D$ introduced
in \cite{Foissy2}, which is a  symmetric and non degenerate Hopf pairing but not a  bidendriform pairing:
$$(\tdun{$d$}\prec\tdun{$d$},\tddeux{$d$}{$d$})=(\tddeux{$d$}{$d$},\tddeux{$d$}{$d$})=0,$$
$$(\tdun{$d$}\otimes \tdun{$d$},\Delta_\prec(\tddeux{$d$}{$d$}))=
(\tdun{$d$}\otimes \tdun{$d$},\tdun{$d$}\otimes \tdun{$d$})=1.$$

The pairing $<,>$ admits a similar description as the pairing  of \cite{Foissy2}.
Let $d\in \D$. We consider the following application:
$$ \xi_d: \left\{\begin{array}{rcl}
\h^\D& \longrightarrow &\h^\D\\
F\in \F^\D & \longrightarrow & 
\left\{\begin{array}{l}
0\mbox{ if the leave of $F$ at most on the right is decorated by $d$;}\\
\mbox{$F$ without its leave at most on the right if it is decorated by $d$.}
\end{array}\right.
\end{array}\right. $$
For example, if $a,b,c,d,e,f\in \D$,
$\xi_f(\tddeux{$a$}{$b$} \tdtroisun{$c$}{$e$}{$d$})=\delta_{e,f}\tddeux{$a$}{$b$}\tddeux{$c$}{$d$}$.

\begin{lemme}
Let $x,y\in \A^\D$, $d\in \D$. Then:
\begin{enumerate} 
\item $<\tdun{$d$}\otimes x, \Delta_\prec(y)>=<x,\xi_d(y)>$.
\item  $<B^+_d(x),y>=<x,\xi_d(y)>$.
\end{enumerate}
\end{lemme}

{\bf Proof.} 

1. We can suppose that $y\in \F^\D$. We denote by $d'$ the decoration of  leave of $y$ at most on the right.
The only eventual admissible cut in ${\cal A}dm_\prec(y)$ 
(defined in proposition \ref{propcoupe}) such that  $P^c(y)=\tdun{$d$}$
is the cut $c_y$ which cut the edge leading to the leave of $y$ at most on the right, if $d=d'$.
Hence, by proposition \ref{propcoupe}:
\begin{eqnarray*}
<\tdun{$d$} \otimes x,y> &=&\sum_{c\in {\cal A}dm_\prec(y)}<\tdun{$d$},P^c(y)><x,R^c(y)>\\
 &=&\sum_{c\in {\cal A}dm_\prec(y)}\delta_{\tdun{$d$},P^c(y)}<x,R^c(y)>\\
 &=&\delta_{d,d'}<x,R^{c_f}(y)>\\
&=&<x,\xi_d(y)>.
\end{eqnarray*}

2. We have
$<B_d^+(x),y>=<\tdun{$d$}\prec x,y>=<\tdun{$d$}\otimes x,\Delta_\prec(y)>=<x,\xi_d(y)>$. $\Box$\\

Hence, the Hopf pairing (also denoted by  $<,>$) induced on $\h^\D$ satisfies the following properties:
\begin{theo}
$<,>:\h^\D\times \h^\D \longrightarrow K$ is the unique pairing such that:
\begin{enumerate}
\item $\forall y\in \h^\D$, $<1,y>=\varepsilon(y)$, where $\varepsilon$ is the counit of $\h^\D$.
\item $\forall x,y\in \h^\D$, $<B_d^+(x),y>=<x,\xi_d(y)>$.
\item $\forall x,y,z\in \h^\D$, $<xy,z>=<x\otimes y,\Delta(z)>$.
\end{enumerate}
Moreover, $<,>$ is symmetric, homogenenous and non-degenerate, and for all $x,y\in \h^\D$, we have
$<S(x),y>=<x,S(y)>$.
\end{theo}

{\bf Remark.} The three points of this theorem allow to compute $<F,G>$ for $F,G \in \F^\D$,
by induction on $weight(F)$.\\

We now give a combinatorial interpretation of this pairing, inspired by our work in \cite{Foissy2}.
If $F \in \F^\D$,
We have two partial orders $\geq_{high}$ et $\geq_{left}$
on the set  $vert(F)$ of vertices of $F$ (see \cite{Foissy2} for more details).
if $x,y$ are two vertices of $F$, then:
\begin{enumerate}
\item We will say that $x \geq_{high} y$ if there is a path from $y$ to $x$.
\item If $s$ and $s'$ are not comparable for $\geq_{high}$, then one of this two vertices 
(say for example $x$) is more on the left than the other: 
we will denote this situation by $x\geq_{left} y$.
\end{enumerate}
These partial orders induce a total order $\geq_{d,l}$ on $vert(F)$ defined by:
$$x\geq_{d,l} y \mbox{ if }(x\geq_{left} y)\mbox{ or } (y\geq_{high} x).$$

\begin{theo}
Let $F,G \in \F^\D$. 
Let ${\cal I}(F,G)$ be the set of bijections $f$ from $vert(F)$ to $vert(G)$ such that: 
\begin{enumerate}
\item $ \forall x,y \in vert(F)$, $x\geq_{high} y$ $\Rightarrow$ $f(x) \geq_{d,l} f(y),$
\item  $\forall x,y \in vert(F)$,$ f(x)\geq_{high} f(y)$ $\Rightarrow$ $x \geq_{d,l} y,$
\item $\forall x \in vert(F)$, $x$ and $f(x)$ have the same decoration.
\end{enumerate}
Then $<F,G>=card({\cal I}(F,G))$.
\end{theo}

{\bf Proof.} 
This is true if $weight(F)\neq weight(G)$, as then $<F,G>=0$ 
and there is no bijection from $vert(F)$ to $vert(G)$.
We suppose that $weight(F)=weight(G)=n$, and let us proceed by induction on $n$. If $n=1$, 
then $<\tdun{$d$},\tdun{$d'$}>=\delta_{d,d'}=card({\cal I}(\tdun{$d$},\tdun{$d'$}))$ by condition 3. 
Suppose that the property is satisfied for all $k<n$
and let $F,G\in \F^\D$ with $n$ vertices. We put $F=t_1 \ldots t_m$.

If $m=1$: then $F=t_1=B_d^+(F')$. So $<F,G>=<F',\xi_d(G)>$.
Let $f \in {\cal I}(F,G)$. Denote by $r$ the root of $F$.
Remark that for all $x \in vert(F)$, $x \geq_{high} r$. By condition 1, 
$f(r)$ is the smallest element of $vert(G)$ for$\geq_{d,l}$:
it is the leave of $G$ at most on the right. If it is not decorated by $d$,
then $<F,G>=0$ and $I(F,G)=\emptyset$ by condition $3$. 
Otherwise, there is a bijection 
${\cal I}(F,G) \longrightarrow {\cal I}(F',\xi_d(G))$, sending $f$ to its restriction \`a $vert(F')$. 
As $<F,G>=<F',\xi_d(G)>$, we have the result.

If $m>1$: we put $F'=t_1 \ldots t_{m-1}$. Let $f\in {\cal I}(F,G)$; consider $f(som(t_m))$. 
Let $y_2$ be a vertex of $G$, such that there exists $x_1\in vert(t_m)$, $f(x_1)\geq_{high} y_2$.
There exists a unique $x_2\in som(F)$, such that $f(x_2)=y_2$.
By condition $2$, $x_1\geq_{d,l} x_2$ and then $x_2\in som(t_m)$. 
Hence, there exists a unique admissible cut $c_f$ of $G$ such that
$R^{c_f}(G)=f(t_m)$ and $P^{c_f}(F)=f(F')$.
We then have a bijection:
$$  \left\{\begin{array}{rcl}
{\cal I}(F,G)&\longrightarrow &\displaystyle
\bigcup_{c \in {\cal A}dm(G)} {\cal I}(F',P^c(G))\times {\cal I}(t_m,R^c(G))\\
f&\longrightarrow&( f_{\mid som(F')} , f_{\mid som(t_m)}) 
\in  {\cal I}(F',P^{c_f}(G))\times {\cal I}(t_m,R^{c_f}(G)).
\end{array}\right. $$

Hence:
\begin{eqnarray*}
<F,G>&=&\sum_{c \in {\cal A}dm(G)} <F',P^c(G)><t_m,R^c(G)>\\
&=&\sum_{c \in {\cal A}dm(G)} <F',P^c(G)><t_m,R^c(G)>\\
&=&\sum_{c \in {\cal A}dm(G)} card({\cal I}(F',P^c(G)))card({\cal I}(t_m,R^c(G)))\\
&=&card\left(\bigcup_{c \in {\cal A}dm(G)} {\cal I}(F',P^c(G))\times {\cal I}(t_m,R^c(G))\right)\\
&=&card({\cal I}(F,G)). \:\Box\\
\end{eqnarray*}

We consider now the case where $\D$ is reduced to a single element. We can identify $\h^\D$ and $\h$.
We will denote $B^+$ instead of $B^+_d$, $\xi$ instead of $\xi_d$.
We define $l_k\in \F$ by induction on $k$ in the following way: $l_0=1$, $l_{k+1}=B^+(l_k)$.
For example:
$$l_1=\tun,\:l_2=\tdeux,\:l_3=\ttroisdeux,\:l_4=\tquatrecinq,\:l_5=\tcinqquatorze\ldots$$

\begin{prop}
Let $F\in \F$, of weight $n$. Then:
\begin{enumerate}
\item $<F,l_n>=1$. 
\item $<F,l_{n-1}\tun>$ is the number of roots of $F$.
\item $<F,\tun l_{n-1}>$ is the number of leaves of $F$.
\end{enumerate}
\end{prop}

{\bf Proof.}

1. Induction on  $n$. If $n=0$, then $F=1$ and $<F,l_0>=<1,1>=1$.
Suppose that the result is true for all forest $G$ of weight $n-1$.
Then $<F,l_n>=<\xi(F),l_{n-1}>=1$ as $\xi(F)\in \F$.\\

2. Let $X$ be the set of $F\in \F$, $F\neq 1$, satisfying $2$.
Let $G\in \F$, with $n$ vertices. We have 
$<B^+(G),l_n\tun>=<G,\xi(l_n\tun)>=<G,l_n>=1$,
so $B^+(G)\in X$: $X$ contents the trees.

Let $F_1,F_2\in X$, with respectively $n_1$ and $n_2$ vertices. We put $n=n_1+n_2$. We have:
\begin{eqnarray*}
<F_1F_2,l_{n-1}\tun>&=&<F_1\otimes F_2,\Delta(l_{n-1}\tun)>\\
&=&\sum_{i+j=n-1} <F_1\otimes F_2, l_i\tun \otimes l_j+l_i \otimes l_j\tun>\\
&=&<F_1\otimes F_2,l_{n_1-1}\tun \otimes l_{n_2}+l_{n_1}\otimes l_{n_2-1}\tun>\mbox{ by homogeneity},\\
&=& <F_1,l_{n_1-1}\tun>+ <F_2,l_{n_2-1}\tun>.
\end{eqnarray*}
So $F_1F_2\in X$. As $X$ contents the trees, $X=\F-\{1\}$.\\

3. Let $Y$ be the set of $F\in \F$, $F\neq 1$, satisfying $3$.
We have $<\tun,\tun l_0>=<\tun,\tun>=1$, so $\tun \in Y$.

Let $G\in Y$, with $n$ vertices. We have 
$<B^+(G),\tun l_n>=<G,\xi(\tun l_n)>=<G,\tun l_{n-1}>$.
As $G$ and $B^+(G)$ have the same number of leaves, $B^+(G)\in Y$, so $Y$ is stable under $B^+$.

Let $F_1,F_2\in Y$, with respectively $n_1$ and $n_2$ vertices. We put $n=n_1+n_2$. We have:
\begin{eqnarray*}
<F_1F_2,\tun l_{n-1}>&=&<F_1\otimes F_2,\Delta(\tun l_{n-1})>\\
&=&\sum_{i+j=n-1} <F_1\otimes F_2, \tun l_i\otimes l_j+l_i \otimes \tun l_j>\\
&=&<F_1\otimes F_2,\tun l_{n_1-1}\otimes l_{n_2}+l_{n_1}\otimes \tun l_{n_2-1}>\mbox{ by homogeneity},\\
&=& <F_1,\tun l_{n_1-1}>+ <F_2,\tun l_{n_2-1}>.
\end{eqnarray*}
So $F_1F_2\in Y$. As $Y$ is stable under product and under $B^+$
and contents $\tun$, $Y=\F-\{1\}$. $\Box$\\


{\bf Examples.} Here are the values of the pairing $<,>$ of $\h$ on forests of weight smaller than $4$:
$$\begin{array}{c|c}
&\tun\\
\hline
\tun&1
\end{array}
\hspace{1cm}
\begin{array}{c|cc}
&\tun\tun&\tdeux\\
\hline
\tun\tun &2&1\\
\tdeux&1&1
\end{array}
\hspace{1cm}
\begin{array}{c|ccccc}
&\tun\tun\tun&\tun\tdeux&\tdeux\tun&\ttroisun&\ttroisdeux\\
\hline
\tun\tun\tun&6&3&3&2&1\\
\tun\tdeux&3&2&2&2&1\\
\tdeux\tun&3&2&2&1&1\\
\ttroisun&2&2&1&1&1\\
\ttroisdeux&1&1&1&1&1
\end{array}$$

$$\begin{array}{c|cccccccccccccc}
&\tun\tun\tun\tun&\tdeux\tun\tun&\tun\tdeux\tun&\ttroisun\tun&
\ttroisdeux\tun&\tun\tun\tdeux&\tdeux\tdeux&\tun\ttroisun&\tun\ttroisdeux&
\tquatreun&\tquatredeux&\tquatretrois&\tquatrequatre&\tquatrecinq\\
\hline
\tun\tun\tun\tun&24&12&12&8&4&12&6&8&4&6&3&3&2&1\\
\tdeux\tun\tun&12&7&7&4&3&7&4&4&3&3&2&2&1&1\\
\tun\tdeux\tun&12&7&7&5&3&7&4&5&3&3&2&2&2&1\\
\ttroisun\tun&8&4&5&3&2&6&3&3&3&2&1&2&1&1\\
\ttroisdeux\tun&4&3&3&2&2&3&2&2&2&1&1&1&1&1\\
\tun\tun\tdeux&12&7&7&6&3&7&4&6&3&6&3&3&2&1\\
\tdeux\tdeux&6&4&4&3&2&4&3&3&2&3&2&2&1&1\\
\tun\ttroisun&8&4&5&3&2&6&3&5&3&3&2&2&2&1\\
\tun\ttroisdeux&4&3&3&3&2&3&2&3&2&3&2&2&2&1\\
\tquatreun&6&3&3&2&1&6&3&3&3&2&1&2&1&1\\
\tquatredeux&3&2&2&1&1&3&2&2&2&1&1&1&1&1\\
\tquatretrois&3&2&2&2&1&3&2&2&2&2&1&2&1&1\\
\tquatrequatre&2&1&2&1&1&2&1&2&2&1&1&1&1&1\\
\tquatrecinq&1&1&1&1&1&1&1&1&1&1&1&1&1&1
\end{array}$$

\bibliographystyle{amsplain}
\bibliography{isomorphism}

\providecommand{\bysame}{\leavevmode\hbox to3em{\hrulefill}\thinspace}
\providecommand{\MR}{\relax\ifhmode\unskip\space\fi MR }
\providecommand{\MRhref}[2]{%
  \href{http://www.ams.org/mathscinet-getitem?mr=#1}{#2}
}
\providecommand{\href}[2]{#2}
\begin{thebibliography}{10}

\bibitem{Abe}
E.~Abe, \emph{Hopf algebras}, Cambridge Tracts in Mathematics, 74, Cambridge
  University Press, Cambridge-New York, 1980.

\bibitem{Aguiar}
M.~Aguiar, \emph{Infinitesimal bialgebras, pre-lie and dendriform algebras},
  Lecture Notes in Pure and Appl. Math., no. 237, Dekker, New York, 2004,
  math.QA/02 11074.

\bibitem{Aguiar2}
Marcelo Aguiar and Frank Sottile, \emph{Structure of the
  {M}alvenuto-{R}eutenauer {H}opf algebra of permutations}, Adv. Math.
  \textbf{191} (2005), no.~2, 225--275, math.CO/02 03282.

\bibitem{Connes}
A.~Connes and D.~Kreimer, \emph{Hopf algebras, {R}enormalization and
  {N}oncommutative geometry}, Comm. Math. Phys \textbf{199} (1998), no.~1,
  203--242, hep-th/98 08042.

\bibitem{Duchamp}
G.~Duchamp, F.~Hivert, and J.-Y. Thibon, \emph{Some generalizations of
  quasi-symmetric functions and noncommutative symmetric functions}, Springer,
  Berlin, 2000, math.CO/01 05065.

\bibitem{Foissy4}
L.~Foissy, \emph{Les alg\`ebres de {H}opf des arbres enracin\'es d\'ecor\'es},
  Th\`ese de doctorat, Universit\'e de Reims, 2002.

\bibitem{Foissy2}
\bysame, \emph{Les alg\`ebres de {H}opf des arbres enracin\'es, {I}}, Bull.
  Sci. Math. \textbf{126} (2002), 193--239.

\bibitem{Foissy3}
\bysame, \emph{Les alg\`ebres de {H}opf des arbres enracin\'es, {II}}, Bull.
  Sci. Math. \textbf{126} (2002), 249--288.

\bibitem{Holtkamp}
R.~Holtkamp, \emph{Comparison of {H}opf {A}lgebras on {T}rees}, Arch. Math.
  (Basel) \textbf{80} (2003), no.~4, 368--383.

\bibitem{Joseph}
A.~Joseph, \emph{Quantum groups and their primitive ideals}, Springer-Verlag,
  1995.

\bibitem{Kreimer1}
D.~Kreimer, \emph{On the {H}opf algebra structure of pertubative quantum field
  theories}, Adv. Theor. Math. Phys. \textbf{2} (1998), no.~2, 303--334,
  q-alg/97 07029.

\bibitem{Kreimer2}
\bysame, \emph{On {O}verlapping {D}ivergences}, Comm. Math. Phys. \textbf{204}
  (1999), no.~3, 669--689, hep-th/98 10022.

\bibitem{Kreimer3}
\bysame, \emph{Combinatorics of (pertubative) {Q}uantum {F}ield {T}heory},
  Phys. Rep. \textbf{4--6} (2002), 387--424, hep-th/00 10059.

\bibitem{Loday4}
J.-L. Loday, \emph{Generalized bialgebras and triples of operads}, avalaible at
  {http://www-irma.u-strasbg.fr/$\tilde{\hspace{.3cm}}$loday/}.

\bibitem{Loday2}
\bysame, \emph{Dialgebras}, Lecture Notes in Math., no. 1763, Springer, Berlin,
  2001, math.QA/01 02053.

\bibitem{Loday}
J.-L. Loday and M.~Ronco, \emph{Hopf algebra of the planar binary trees}, Adv.
  Math. \textbf{139} (1998), no.~2, 293--309.

\bibitem{Loday3}
\bysame, \emph{Order structure on the algebra of permutations and of planar
  binary trees}, J. Algebraic Combin. \textbf{15} (2002), no.~3, 253--270,
  math.CO/01 02066.

\bibitem{Malvenuto}
C.~Malvenuto and Ch. Reutenauer, \emph{Duality between quasi-symmetric
  functions and the {S}olomon descent algebra}, J. Algebra \textbf{177} (1995),
  no.~3, 967--982.

\bibitem{Markl}
M.~Markl, S.~Shnider, and J.~Stasheff, \emph{Operads in algebra, topology and
  physics}, Mathematical Surveys and Monographs, 96, American Mathematical
  Society, Providence, RI, 2002.

\bibitem{Milnor}
J.W. Milnor and J.~C. Moore, \emph{On the structure of {H}opf algebras}, Ann.
  of Math. (2) \textbf{81} (1965), 211--264.

\bibitem{Ronco}
M.~Ronco, \emph{Primitive elements of a free dendriform algebra}, Contemp.
  Math. \textbf{267} (2000), 245--263.

\bibitem{Ronco2}
\bysame, \emph{Eulerian idempotents and {M}ilnor-{M}oore theorem for certain
  non-cocommutative {H}opf algebras}, J. Algebra \textbf{254} (2002), no.~1,
  152--172.

\end{thebibliography}

\end{document}